\documentclass[11pt]{article}
\usepackage{latexsym,amsfonts,amssymb,amsmath,euscript,amsthm}
\usepackage[applemac]{inputenc}  
\usepackage[english,francais]{babel}  
\usepackage{graphicx} 
\usepackage{subfigure}
\usepackage{hyperref}         
\usepackage{color,fancybox}
\usepackage{bm}
\usepackage{float}

\numberwithin{equation}{section}
\newtheorem{teo}{Theorem}[section]
\newtheorem{lem}[teo]{Lemma}
\newtheorem{cor}[teo]{Corolary}
\newtheorem{prop}[teo]{Proposition}
\newtheorem{defi}[teo]{Definition}
\newtheorem{re}[teo]{Remark}

\newcommand{\A}{\mathcal{A}_0}

\newcommand{\eps}{\varepsilon}

\newcommand{\eto}{\stackrel{\eps\to 0}{\longrightarrow}}

\newcommand{\R} {\mathbb{R}}
\newcommand{\N} {\mathbb{N}}

\def\eto{\buildrel \epsilon\to 0\over\longrightarrow }

\title{Distance of attractors for thin domains\footnote{
This research has been partially supported by grants MTM2016-75465,  MTM2012-31298, ICMAT Severo Ochoa project SEV-2015-0554 (MINECO), Spain and Grupo de Investigaci\'on CADEDIF, UCM.}   }
\author{Jos\'{e} M. Arrieta\footnote{Departamento de Matem\'atica Aplicada, Universidad Complutense de Madrid, 28040 Madrid and Instituto de Ciencias Matem\'aticas
CSIC-UAM-UC3M-UCM, C/Nicol\'as Cabrera 13-15, Cantoblanco, 28049 Madrid, Spain}
and Esperanza Santamar\'ia\footnote{Universidad a Distancia de Madrid y Colegio San Patricio, Madrid, Spain} 
}

\date{ }

\begin{document}

\maketitle
%
%
%
%

{\footnotesize 
\par\noindent {\bf Abstract:}
In this work we consider a dissipative reaction-diffusion equation in a $d$-dimensional thin domain shrinking to a one dimensional segment and obtain good rates for the convergence of the attractors. To accomplish this, we use estimates on the convergence of inertial manifolds as developed previously in \cite{Arrieta-Santamaria-C0} and Shadowing theory.  

 \vskip 0.5\baselineskip

\vspace{11pt}

\vspace{6pt}

}

\numberwithin{equation}{section}
\newtheorem{theorem}{Theorem}[section]
\newtheorem{lemma}[theorem]{Lemma}
\newtheorem{corollary}[theorem]{Corollary}
\newtheorem{proposition}[theorem]{Proposition}
\newtheorem{definition}[theorem]{Definition}
\newtheorem{remark}[theorem]{Remark}
\allowdisplaybreaks

\section{Introduction}
\selectlanguage{english}

In this work we study the rate of convergence of attractors for a reaction diffusion equation in a thin domain when the thickness of the domain goes to zero. Our domain is a thin channel obtained by shrinking a fixed domain $Q\subset \mathbb{R}^d$, see Figure \ref{dominioQ}, by a factor $\eps$ in $(d-1)$-directions. The thin channel $Q_\eps$ collapses to the one dimensional line segment $[0, 1]$ as $\eps$ goes to zero. 

The reaction diffusion-equation in $Q_\eps$ is given by 
\begin{equation}\label{equationonQepsilon}
\left\{
\begin{array}{r@{=} l c}
u_t-\Delta u+ \mu u \;&\;f(u)\quad &\textrm{in}\quad Q_\varepsilon,\\
\frac{\partial u}{\partial\nu_\varepsilon} \;&\;0\quad&\textrm{in}\quad \partial Q_\varepsilon,
\end{array}
\right.
\end{equation}
where $\mu>0$ is a fixed number, $\nu_\varepsilon$ the unit outward normal to $\partial Q_\varepsilon$ and $f:\mathbb{R}\rightarrow\mathbb{R}$ is a nonlinear term, with appropriate dissipativity conditions to guarantee the existence of an attractor $\mathcal{A}_\eps\subset H^1(Q_\eps)$.   

As the parameter $\eps\to 0$, the thin domain shrinks to the line segment $[0,1]$ and the limiting reaction-diffusion equation is given by
\begin{equation}\label{equationon(01)}
\left\{
\begin{array}{r@{=} l c}
u_t-\frac{1}{g}(gu_x)_x+\mu u\;&\;f(u)\qquad&\textrm{in}\quad (0, 1),\\
u_x(0)=u_x(1)\;&\;0.
\end{array}
\right.
\end{equation}
which also has an attractor $\mathcal{A}_0\subset H^1(0,1)$. 

There are several works in the literature comparing the dynamics of both equations 
and showing the convergence of $\mathcal{A}_\eps$ to $\mathcal{A}_0$ as $\eps\to 0$, under certain hypotheses. One of the most relevant and pioneer work in this direction is \cite{Hale&Raugel3}, where the authors show that when $d=2$ and every equilibrium of the limit problem \eqref{equationon(01)} is hyperbolic, the attractors behave continuously and moreover, the flow in the attractors of both systems are topologically conjugate.  In order to accomplish this task, the authors exploit the fact that the limit problem is one dimensional, which allows them to construct inertial manifolds for \eqref{equationonQepsilon} and \eqref{equationon(01)} which will be close in the $C^1$ topology. Restricting the flow to these inertial manifolds, and using that the limit problem is Morse-Smale (under the condition that all equilibria being hyperbolic, see \cite{Henry2}) they prove the $C^0$-conjugacy of the flows. Moreover the method of constructing the inertial manifolds for fixed $\eps\in [0,\eps_0]$ consists in using the method described in \cite{Mallet-Sell}. They consider the finite dimensional linear manifold given by the span of the eigenfunctions corresponding to the first $m$ eigenvalues of the elliptic operator and let evolve this linear manifold with the nonlinear flow, which $\omega$-limit set is a $C^1$ manifold and it is the inertial manifold,  which, as a matter of fact it is a graph over the finite dimensional linear manifold. This method provides them with an estimate of the distance of the inertial manifolds of the order of $\eps^\gamma$ for some $\gamma<1$.   Later on, reducing the system to the inertial manifolds and using the general techniques to estimate the distance of attractors for gradient flows, see \cite{Hale&Raugel} Theorem 2.5, give them the estimate $\eps^{\gamma'}$ with some $\gamma'<\gamma<1$ which depends on the number of equilibria of the limit problem and other characteristics of the problem.

Our setting is more general than the one from \cite{Hale&Raugel3}, since we consider general $d-$dimensional thin domains (not just 2-dimensional).   Moreover, our approach to this problem has some differences with respect to theirs.  In our case, we will also construct inertial manifolds, but we will construct them following the Lyapunov-Perron method, as developed in \cite{Sell&You}. This method, as it is shown in \cite{Arrieta-Santamaria-C0,Arrieta-Santamaria-C1} provides us with good estimate of the $C^0$ distance of the inertial manifolds (which is of order $\eps|\ln(\eps)|$, see \cite{Arrieta-Santamaria-C0})  and  with the $C^{1,\theta}$ convergence of this manifolds, see \cite{Arrieta-Santamaria-C1}.

%
%

Once the Inertial Manifolds are constructed and we have a good estimate of its distance we can project the systems to these inertial manifolds and obtain the reduced systems, which are finite dimensional. The limit reduced system will be a Morse-Smale gradient like system, see \cite{Hale}. Then  Shadowing theory and its relation to the distance of the attractors, as developed in Appendix \ref{shadowing}, will give us the key to obtain the rates of convergence of the attractors.   

Let us mention that the estimate we find on the Hausdorff symmetric  distance of the attractors is the following (see Theorem \ref{maintheorem-thindomain}), 
$$\hbox{dist}_{H^1(Q_\eps)} (\mathcal{A}_0, \mathcal{A}_\eps)\leq  C \eps^{\frac{d+1}{2}}|\log(\eps)|$$
which improves the one obtained in \cite{Hale&Raugel3}.

We describe now the contents of this chapter:

In section  \ref{setting} we give a complete description of the thin domain $Q_\eps$, will set up the basic notation we will need. We also introduce the main result of the paper.

In section \ref{elliptic} we study the related elliptic problem, obtaining an estimate for the distance of the resolvent operators and proving this estimate is optimal.  We postpone the proof of the main result of this section Proposition \ref{resolvente} to Appendix \ref{proof-resolvente}.  

In section \ref{nonlinear} we analyze the nonlinearity and we prepare it for the construction of inertial manifolds. We make an appropriate cut off of the non-linear term and analyze the conditions this new nonlinearity satisfies. 

In section \ref{InertialManifoldsConstruction} we construct the corresponding inertial manifolds, reducing our problem to a finite dimensional one.

In section \ref{distanceattractors} using the estimates on the distance of the inertial manifolds together with the shadowing result  obtained in Appendix \ref{shadowing} we provide an almost optimal rate of convergence of attractors, proving the main theorem Theorem \ref{maintheorem-thindomain}. 

At the end we have included two appendices. Appendix \ref{proof-resolvente} contains the proof of Proposition \ref{resolvente} and Appendix \ref{shadowing} contains some results on the relation of Shadowing and the distance of attractors for Morse-Smale maps.

\section{Setting of the problem and main results}
\label{setting}
In this section we set up the problem, describing clearly the domain and the equations we are dealing with. We will also state our main result on the distance of attractors. We end up the section with some notation and technical results needed thereafter.

We start describing the thin domain. Let $\Omega= (0, 1)$ and let $Q$ be the set
$$Q=\{(x, \mathbf{y})\in\mathbb{R}^d: 0\leq x\leq1,\; \; \mathbf{y}\in\Gamma^1_x\},$$
with $d\geq2$, and $\Gamma_x^1$ diffeomorphic to the unit ball in $\mathbb{R}^{d-1}$, $B(0, 1)$, for all $x\in[0,1]$, see Figure \ref{dominioQ}, that is, we assume that for each $x\in[0, 1]$, there exists a $C^1$ dipheomorphism $\mathbf{L}_x$
\begin{equation}\label{difeomorfismo}
\mathbf{L}_x:B(0,1)\longrightarrow\Gamma_x^1\subset \mathbb{R}^{d-1}.
\end{equation}
We also assume that, if we define
\begin{equation}\label{definition-of-J}
\left\{
\begin{array}{r c l }
\mathbf{L}:(0,1)\times B(0, 1)\; &\;\longrightarrow \;&\;Q\\
(x,\mathbf{y})\;&\;\mapsto\;&\;(x, \mathbf{L}_x(\mathbf{y}))
\end{array}
\right.
\end{equation}
then $\mathbf{L}$ is a $C^1$ diffeomorphism. The boundary of $Q$ has two distinguished parts, the one formed by $\Gamma_0^1\cup \Gamma_1^1$ (the two lids of the thin domain) and the lateral boundary $\partial_lQ=\{ (x,y):  x\in (0,1), y\in \partial \Gamma_x^1\}$

\begin{figure}\label{dominioQ}
\centerline{\includegraphics[scale=0.6]{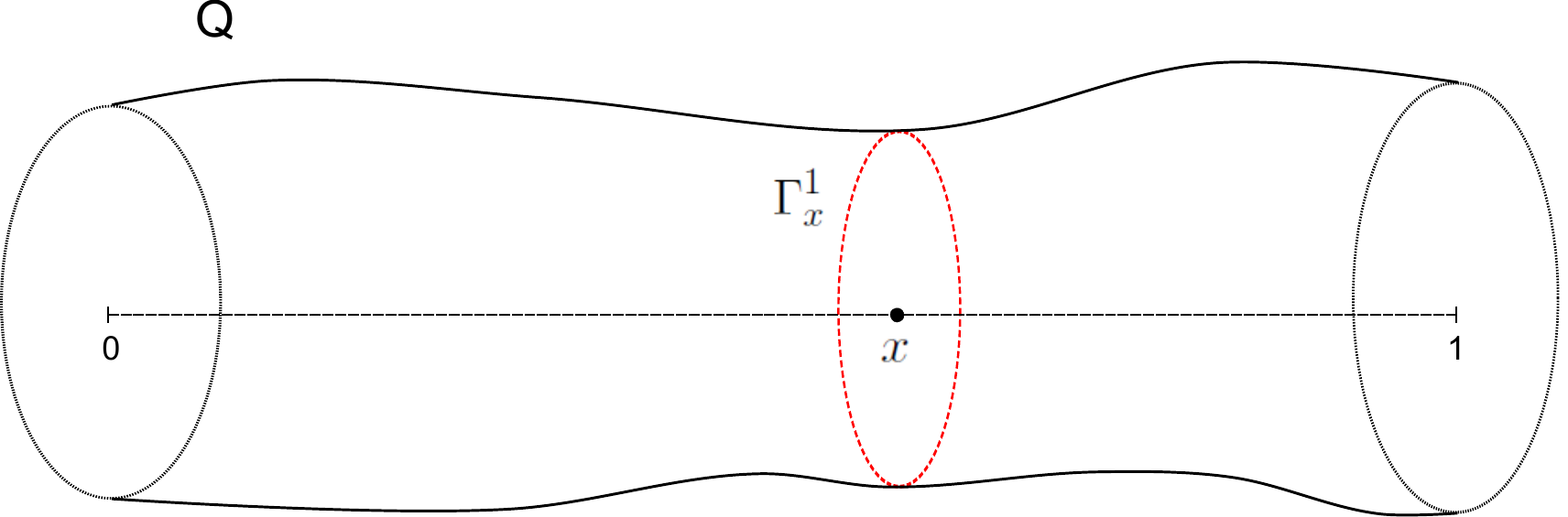}}
\caption{ Domain $Q$ with $d=3$.\label{dominioQ}}
\end{figure}

 Our thin channel, or thin domain will be defined by 
 $$Q_\varepsilon =\{(x, \varepsilon\mathbf{y})\in\mathbb{R}^d: (x, \mathbf{y})\in Q\},\qquad \varepsilon\in (0,1).$$

Notice that this set is obtained by shrinking the set $Q$ by a factor $\varepsilon$ in the $(d-1)$-directions given by the variable $\mathbf{y}\in \mathbb{R}^{d-1}$.  This domain gets thinner and thinner as $\varepsilon\to 0$ and it approaches the one dimensional line segment given by $\Omega\times\{\mathbf{0}\}=(0,1)\times\{\mathbf{0}\}$. 

\par\medskip 

We denote by $g(x):=|\Gamma_x^1|$ the $(d-1)$-dimensional Lebesgue measure of the set $\Gamma_x^1$. From the hypothesis of the smoothness of the map $\mathbf{L}$ above, see \eqref{definition-of-J},  we have that $g$ is a smooth function defined in $[0, 1]$. In particular, there exist $g_0, g_1>0$ such that $g_0\leq g(x)\leq g_1$ for all $x\in [0, 1]$. 

\begin{re}
An important subclass of these thin domains are those whose transversal sections $\Gamma_x^1$ are disks centered at the origin of radius $r(x)$, that is, 
$$Q=\{(x, \mathbf{y})\in\mathbb{R}^d: 0\leq x\leq 1, |\mathbf{y}|< r(x)\}.$$
In this particular case, $g(x)=|B(0,1)|r(x)^{d-1}$, with $|B(0,1)|$ the Lebesgue measure of the unit ball in $\mathbb{R}^{d-1}$. The diffeomorphism $\mathbf{L}$ defined in (\ref{definition-of-J}) is given by,
$$\mathbf{L}(x, \mathbf{y})=(x, r(x)\mathbf{y}).$$
\end{re}

\vspace{0.5cm}

We consider the following reaction-diffusion equation in $Q_\varepsilon$, $0<\varepsilon\leq\varepsilon_0$,
\begin{equation}\label{equationonQepsilon}
\left\{
\begin{array}{r@{=} l c}
u_t-\Delta u+ \mu u \;&\;f(u)\quad &\textrm{in}\quad Q_\varepsilon,\\
\frac{\partial u}{\partial\nu_\varepsilon} \;&\;0\quad&\textrm{in}\quad \partial Q_\varepsilon,
\end{array}
\right.
\end{equation}
where $\mu>0$ is a fixed number, $\nu_\varepsilon$ the unit outward normal to $\partial Q_\varepsilon$ and $f:\mathbb{R}\rightarrow\mathbb{R}$ a $C^2$-function satisfying the following  growth condition 
\begin{equation}\label{growthcondition}
|f'(s)|\leq C(1+|s|^{\rho-1}), \quad s\in \mathbb{R}
\end{equation}
for some $\rho\geq 1$, 
and the dissipative condition,
\begin{equation}\label{dissipativecondition1}
\exists\; M>0,\quad\textrm{s.t.}\quad f(s)\cdot s\leq 0,\qquad |s|\geq M.
\end{equation}

With the growth condition \eqref{growthcondition} we know that problem \eqref{equationonQepsilon} is locally well posed in some functional space of the type $L^r(Q_\eps)$ for some $r>1$, maybe large enough, or $W^{1,p}(Q_\eps)$, see \cite{Arrieta+CarvalhoTAMS, Arrieta+Carvalho+Anibal}.  With the dissipative condition  and with some regularity arguments we obtain that solutions are globally defined and with the aid of the maximum principle there exist uniform asymptotic bounds in the
sup norm of the solutions. That is, for any initial condition $\phi_\eps$ there exists a time $\tau$, that may depend on $\eps$ and on the initial condition, such that the solution starting at $\phi_\eps$ after time $\tau$ is uniformly bounded by $M$, that is $|u(t,x,\phi_\eps)|\leq M$ for $t\geq \tau$, with $M$ given by \eqref{dissipativecondition1}. This uniform asymptotic bounds together with parabolic regularity theory imply that the equation \eqref{equationonQepsilon} has an attractor $\mathcal{A}_\eps\subset H^1(Q_\eps)\cap L^\infty(Q_\eps)$ satisfying the uniform bound 
\begin{equation}\label{uniformbound-eps}
\|u_\eps\|_{L^\infty(Q_\eps)}\leq M, \qquad \hbox{for all }u_\eps \in \mathcal{A}_\eps
\end{equation}

The limit problem of (\ref{equationonQepsilon}) is given by, see \cite{Hale&Raugel3},
\begin{equation}\label{equationon(01)}
\left\{
\begin{array}{r@{=} l c}
u_t-\frac{1}{g}(gu_x)_x+\mu u\;&\;f(u)\qquad&\textrm{in}\quad (0, 1),\\
u_x(0)=u_x(1)\;&\;0.
\end{array}
\right.
\end{equation}
and, just as the analysis above, this equation has also an attractor $\mathcal{A}_0\subset H^1(0,1)\cap L^\infty(0,1)$ satisfying also the bounds 
\begin{equation}\label{uniformbound0}
\|u_0\|_{L^\infty(0,1)}\leq M.
\end{equation}
Observe that the dynamical system generated by this equation has a gradient structure (see \cite{Hale}) and in particular its attractor is formed by equilibria and connections among them. Moreover, if all equilibria are hyperbolic then we have only a finite number of them and the system has a Morse-Smale structure (see \cite{Henry2}). 

 Notice that in a natural way we may consider the attractor $\mathcal{A}_0$ as a subset of $H^1(Q_\eps)$, just by considering that any function $u_0(x)$ defined in $(0,1)$ is extended to all of $Q_\eps$ by 
$\tilde u_0(x,\mathbf{y})=u_0(x)$.

\par\bigskip

We now introduce the main result of the paper.
\begin{teo}\label{maintheorem-thindomain}
Under the notations above and assuming that all equilibria of problem \eqref{equationon(01)} are hyperbolic, we have
\begin{equation}\label{distance-attractors}
\hbox{dist}_{H^1(Q_\eps)} (\mathcal{A}_0, \mathcal{A}_\eps)\leq  C \eps^{\frac{d+1}{2}}|\log(\eps)|,
\end{equation}
with $\hbox{dist}_X(\cdot,\cdot)$ the symmetric Haussdorf distance in the space $X$. 
\end{teo}

Recall that $\hbox{dist}_X(A,B)$ is defined as 
$$\hbox{dist}_X(A,B)=\max\{ \sup_{x\in A} \inf_{y\in B} d(x,y), \sup_{y\in B} \inf_{x\in A} d(x,y)\}.$$


%
%
%
%

\par\bigskip

Next, we present the notation and some conditions needed for the proof. 

As we have noted above, the attractors of both equations \eqref{equationonQepsilon} and  \eqref{equationon(01)} have
uniform $L^\infty$ bounds, as expressed in \eqref{uniformbound-eps} and 
\eqref{uniformbound0}. This fact will allow us to cut off the nonlinearity $f$ outside the interval $(-M,M)$ so that the new nonlinearity that we will still denote by $f$ has compact support and coincides with the old one in $(-M,M)$, satisfies  
\begin{equation}\label{dissipativecondition2}
|f(s)|+|f'(s)|+|f''(s)|\leq L_f \qquad\textrm{for all}\quad s\in\mathbb{R},
\end{equation}
and the dissipative condition \eqref{dissipativecondition1} still holds for the new $f$. Moreover, since the attractors for the old nonlinearity satisfy \eqref{uniformbound-eps} and \eqref{uniformbound0} and the new $f$ coincides with the old one in $(-M,M)$, then the attractors for the new equations are exactly the same as the attractors for the original equations. This means that we may assume from the beginning that the nonlinearity $f$ satisfies \eqref{dissipativecondition2}

When dealing with problems where the domain varies it is sometimes convenient to make transformations, as simple as possible, so that we 
transform all problems to a fixed reference domain.  This will imply in many instances that the parameter appears in the equation and usually it will 
show up as a singular parameter.  In our case, we will transform problem 
\eqref{equationonQepsilon} into a problem in the fixed set $Q=\{(x, \mathbf{y})\in\mathbb{R}^d: 0\leq x\leq1, \mathbf{y}\in\Gamma_x^1\}$, (Figure \ref{dominioQ}).  The transformation we will use is $(x,\mathbf{y})\to (x,\frac{\mathbf{y}}{\eps})$.
 With this transformation, the reaction-diffusion equation (\ref{equationonQepsilon}) is transformed into the following equation on the fixed domain $Q$,
\begin{equation}\label{equationonQ}
\left\{
\begin{array}{r@{=} l c}
u_t-\frac{\partial^2u}{\partial x^2}-\frac{1}{\varepsilon^2}\Delta_{\mathbf{y}} u+\mu u\;&\;f(u)\qquad&\textrm{in}\quad Q,\\
\frac{\partial u}{\partial\nu_x}+\frac{1}{\varepsilon^2}\frac{\partial u}{\partial\nu_\mathbf{y}}\;&\;0\quad&\textrm{on}\quad \partial Q
\end{array}
\right.
\end{equation}
where $\nu=(\frac{\partial u}{\partial\nu_x},\frac{\partial u}{\partial\nu_\mathbf{y}}) $ is the unit outward normal to $\partial Q$.

The natural spaces to analyze (\ref{equationonQ}) are given by,
$$H^1_{\boldsymbol\varepsilon}(Q):=(H^1(Q), \|\cdot\|_{H^1_{\boldsymbol\varepsilon}(Q)}),$$
with the norm
$$\|u\|_{H^1_{\bm\varepsilon}(Q)}:=\left(\int_Q(|\nabla_xu|^2+\frac{1}{\varepsilon^2}|\nabla_{\mathbf{y}}u|^2+|u|^2)dxd\mathbf{y}\right)^{1/2},$$
and $L^2(Q)$ with the usual norm $\|\cdot\|_{L^2(Q)}$.

\par\medskip

Notice that if we define  the isomorphism $\mathbf{i}_{\bm\varepsilon}: L^2(Q_\varepsilon)\rightarrow L^2(Q)$ as 
$$\mathbf{i}_{\bm\varepsilon}(u)(x,\mathbf{y}):= u(x, \varepsilon \mathbf{y}),$$ 
its restriction to $H^1(Q_\varepsilon)$ is also an isomorphism from $H^1(Q_\varepsilon)$ to $H^1(Q)$ (or equivalently to $H^1_{\bm\eps}(Q)$).  Then we easily have the following identities:
\begin{equation}\label{norm-relations1}
\|\mathbf{i}_{\bm\eps}(u)\|_{L^2(Q)}=\eps^{-\frac{d-1}{2}}\|u\|_{L^2(Q_\eps)}
\end{equation}
\begin{equation}\label{norm-relations2}
\|\mathbf{i}_{\bm\eps}(u)\|_{H^1_{\bm\eps}(Q)}=\eps^{-\frac{d-1}{2}}\|u\|_{H^1(Q_\eps)}
\end{equation}

The isomorphism $\mathbf{i}_{\bm\eps}$ also allows us how to relate easily the semigroups generated by \eqref{equationonQepsilon} and \eqref{equationonQ} as follows: if $S_\eps(t)$ is the semigroup generated by \eqref{equationonQepsilon}  and $\tilde S_\eps(t)$ the one from \eqref{equationonQ}, then we have 

%
$$S_\varepsilon(t)(\cdot):=\mathbf{i}^{-1}_{\bm\varepsilon}\circ \tilde{S}_\varepsilon(t)\circ \mathbf{i}_{\bm\varepsilon}(\cdot),$$

The limit problem of eqution (\ref{equationonQ}) is also given by (\ref{equationon(01)}).

The natural spaces to treat the limit problem are the following
$$L^2_g(0,1):=(L^2(0,1), \|\cdot\|_{L^2_g(0,1)})  \qquad\textrm{with}\quad \|u\|_{L^2_g(0,1)}:=\left(\int_0^1g(x)|u(x)|^2dx\right)^{\frac{1}{2}},$$
and
$$H^1_g(0,1):=(H^1(0,1), \|\cdot\|_{H^1_g(0,1)}) \qquad\textrm{with}\quad \|u\|_{H^1_g(0,1)}:= \left(\int_0^1g(x)(|u_x|^2+|u|^2)dx\right)^{\frac{1}{2}}.$$
\vspace{0.2cm}

Throughout this paper we will denote by $|\cdot|$ the norm in $\mathbb{R}^d$.
\vspace{0.3cm}


Both evolution problems (\ref{equationonQ}) and (\ref{equationon(01)}) admit an abstract formulation that we are going to overview here. 
Let $A_\eps:D(A_\eps)\subset L^2(Q)\to L^2(Q),$ with
$$A_\varepsilon=-\frac{\partial^2}{\partial x^2}-\frac{1}{\varepsilon^2}\Delta_y+\mu I,\qquad\textrm{and}\qquad D(A_\eps)=\{u\in H^2(Q): \frac{\partial u}{\partial \nu}=0,\, \textrm{at}\, \partial Q\}$$ 
and
$A_0:D(A_0)\subset L^2_g(0,1)\to L^2_g(0,1),$
with
$$A_0v=-\frac{1}{g}(gv_x)_x+\mu v\qquad\textrm{ and}\qquad D(A_0)=\{ v\in H_g^2(0,1), u'(0)=u'(1)=0\}.$$ 

Both operators are selfadjoint, positive linear operators with compact resolvent and they are defined on separable Hilbert spaces.  We denote by $X_\eps=L^2(Q)$, $0<\eps\leq \eps_0$ with the usual norm and $X_0=L^2_g(0,1)$, with its norm defined above. Let also denote by $X_\eps^1=D(A_\eps)$ with the graph norm and similarly for $X_0^1$.  We also consider the fractional power spaces $X_\eps^\alpha$ for $0<\alpha<1$ and $0\leq \eps<\eps_0$, see \cite{Henry1}. In particular the spaces $X_\eps^{1/2}$ is $H^1_\eps(Q)$ and $X_0^{1/2}$ is $H^1_g(0,1)$ defined above.

Hence, (\ref{equationonQ}) and (\ref{equationon(01)}) can be written as  
\begin{equation}\label{problemaperturbado-thindomain}
(P_\varepsilon)\left\{
\begin{array}{r l }
u^\varepsilon_t+A_\varepsilon u^\varepsilon&=F_\varepsilon(u^\varepsilon),\qquad 0<\varepsilon\leq \eps_0\\
u^\varepsilon(0)\in X^\alpha_\varepsilon,
\end{array}
\right.
\end{equation}
and
\begin{equation}\label{problemalimite-thindomain}
(P_0)\left\{
\begin{array}{r l }
u^0_t+A_0u^0&=F_0(u^0),\\
u^0(0)\in X^\alpha_0,
\end{array}
\right.
\end{equation}
where $F_\eps$ and $F_0$ are the nonlinearity $f$ acting in the appropriate fractional power spaces, which will be analyzed in detail in Section \ref{nonlinear}.

 We define an extension operator which maps functions defined in $[0, 1]$ into functions defined in $Q$. The natural way to construct this operator is to extend the functions defined in $[0,1]$ constantly in the other $d-1$ variables. Therefore we denote by $E$ the transformation,
 \begin{equation}\label{operadorextension}
 \left.
 \begin{array}{r c l}
E: L^2_g(0,1)\;&\;\longrightarrow\;&\;L^2(Q)\\
u\;&\;\mapsto\;&\; E(u) (x,\mathbf{y})=u(x)
 \end{array}
 \right.
 \end{equation}
 In a similar fashion we may define the transformation $E_\varepsilon: L^2_g(0,1)\rightarrow L^2(Q_\varepsilon)$ defined as $(E_\varepsilon u)(x, \mathbf{y})=u(x)$. The difference with $E$ is that $E_\varepsilon$ lands in $L^2(Q_\varepsilon)$. As a mather of fact, $E_\varepsilon= \mathbf{i}_{\bm\varepsilon}^{-1}\circ E$.
 
These transformations can also be considered restricted to $H^1_g(0,1)$. In this case, we have $E:H_g^1(0,1)\longrightarrow H_{\bm\varepsilon}^1(Q)$ and $E_\varepsilon: H^1_g(0, 1)\rightarrow H^1(Q_\varepsilon)$. These transformations can be considered too as $E:X_0^\alpha\longrightarrow X_\eps^\alpha$.

To compare functions from $L^2(Q)$ and $L^2_g(0,1)$ (and from $X_\eps^\alpha$ and $X_0^\alpha$, respectively) we also need a projection operator $M$, defined as follows,
\begin{equation}\label{M-definition}
 \left.
 \begin{array}{r c l}
M: L^2(Q)\;&\;\longrightarrow\;&\;L^2_g(0, 1)\\
u\;&\;\mapsto\;&\; \displaystyle M(u) (x)=\frac{1}{|\Gamma_x^1|}\int_{\Gamma_x^1}u(x, \mathbf{y})d\mathbf{y},
 \end{array}
 \right.
 \end{equation}
similary, we may define the map,
\begin{equation}\label{M-eps-definition}
 \left.
 \begin{array}{r c l}
M_\varepsilon: L^2(Q_\varepsilon)\;&\;\longrightarrow\;&\;L^2_g(0, 1)\\
u\;&\;\mapsto\;&\; \displaystyle  M_\varepsilon(u) (x)=\frac{1}{|\Gamma_x^\varepsilon|}\int_{\Gamma_x^\varepsilon}u(x, \mathbf{y})d\mathbf{y},
 \end{array}
 \right.
 \end{equation}
and, in the same way, for $0<\alpha<\frac{1}{2}$, $M:X_\eps^\alpha\longrightarrow X_0^\alpha$. Moreover $M: H^1_{\bm\eps}(Q)\longrightarrow H^1_g(0,1)$ and $M_\varepsilon: H^1(Q_\varepsilon)\longrightarrow H^1_g(0,1)$.  
\vspace{0.3cm}

The following estimates are straightforward: 

$$\|M\|_{\mathcal{L}(L^2(Q), L^2_g(0,1))}\leq1, \quad \|M\|_{\mathcal{L}(H^1_{\bm\varepsilon}(Q), H^1_g(0,1))}\leq1, \quad \|M_\varepsilon\|_{\mathcal{L}(L^2(Q_\varepsilon), L^2_g(0,1))}\leq \varepsilon^{\frac{1-d}{2}}.$$ 

$$\|E u\|_{L^2(Q)}=\|u\|_{L^2_g(0, 1)},\quad\|E_\varepsilon u\|_{L^2(Q_\varepsilon)}=\varepsilon^{\frac{d-1}{2}}\|u\|_{L^2_g(0, 1)}  \quad \forall u\in L^2_g(0,1)$$ 
$$\|E u\|_{H^1_{\bm\varepsilon}(Q)}=\|u\|_{H^1_g(0, 1)},\quad\forall u\in H^1_g(0,1)$$

\begin{re}\label{EMalpha}

i)  From \cite[Theorem 16.1, pg 528]{Yagi2010} we get that the spaces of fractional power of operators and the spaces obtained via interpolation coincide and even they are isometric.  This means that $X_0^\alpha=[L^2_g(0, 1), H^1_g(0, 1)]_{2\alpha}$ and $X_\eps^\alpha=[L^2(Q), H^1_\eps(Q)]_{2\alpha}$ with isometry.  This implies that we also have 
\begin{equation}\label{Malpha}
\|Mu\|_{X_0^\alpha}\leq \|u\|_{X_\eps^\alpha}  
\end{equation}

%
For the operator $E:X_0^\alpha\longrightarrow X_\eps^\alpha$ we also obtain,
\begin{equation}\label{Ealpha}
\|Eu\|_{X_\eps^\alpha}\leq \|u\|_{X_0^\alpha}, \qquad\forall u\in X_0^\alpha,
\end{equation}
applying exactly the same arguments. 

Note that \eqref{Malpha} and \eqref{Ealpha} show that estimates \eqref{cotaextensionproyeccion} are satisfied with  $\kappa=1$.

\par\bigskip  ii)  Moreover, via interpolation we easily get that if $\alpha<1/2$, we get  $[L^2(Q), H^1_\eps(Q)]_{2\alpha}\hookrightarrow[L^2(Q), H^1(Q)]_{2\alpha}\hookrightarrow L^{\frac{2d}{d-4\alpha}}(Q)$ with an embedding constant independent of $\eps$.  Hence, we also have that the embedding constant of $X_\eps^\alpha\hookrightarrow L^{\frac{2d}{d-4\alpha}}(Q)$ is independent of $\eps$. 

\end{re}
\bigskip

We  include now a technical result on the operator $E$ that will be used later.  
\begin{lem}\label{normaproyeccionextension}
We have the following 
%
\par\noindent i) There exists a constant $\beta>0$ such that 
$$\|u_\varepsilon- E M u_\varepsilon\|^2_{L^2(Q)}\leq \beta \|\nabla_{\mathbf{y}}u_\varepsilon\|^2_{L^2(Q)}, \quad \forall u_\eps\in H^1(Q)$$
$$\|w_\varepsilon-E_\varepsilon M_\varepsilon w_\varepsilon\|^2_{L^2(Q_\varepsilon)}\leq \beta\varepsilon^2\|\nabla_{\mathbf{y}}w_\varepsilon\|^2_{L^2(Q_\varepsilon)}, \quad \forall w_\eps\in H^1(Q_\eps)$$
\par\noindent ii) Let $K\subset X_0^\alpha$ a compact set. Then,
$$\sup_{u_0\in K}\left|\|Eu_0\|_{X_\eps^\alpha}-\|u_0\|_{X_0^\alpha}\right|\rightarrow0,\qquad \textrm{as}\;\; \eps\to 0.$$
\end{lem}
\paragraph{ Proof.}  
%

\par\noindent i)  Observe that,
$$\|u_\varepsilon- E M u_\varepsilon\|_{L^2(Q)}^2=\int_0^1\int_{\Gamma_x^1}|u_\varepsilon(x,\mathbf{y})-(Mu_\varepsilon)(x)|^2d\mathbf{y}dx\leq\int_0^1\frac{1}{\lambda_2(\Gamma_x^1)}\int_{\Gamma_x^1}|\nabla_{\mathbf{y}}u_\varepsilon(x, \mathbf{y})|^2d\mathbf{y}dx$$
where we are using Poincare inequality in $\Gamma_x^1$ ($\lambda_2(\Gamma_x^1)$ is the second Neumann eigenvalue in $\Gamma_x^1$). 
%

Let us see that there exists a $\hat{\lambda}_2>0$ such that,
$$\lambda_2(\Gamma_x^1)\geq\hat{\lambda}_2>0,\qquad\forall x\in [0,1].$$
If this is not the case, then there exists a sequence $x_n\rightarrow x_0\in [0,1]$ such that $\lambda_2(\Gamma_{x_n}^1)\rightarrow 0$ as $n\rightarrow\infty$. But $\Gamma_{x_n}^1$ for $n$ large enough is $C^1$ close to $\Gamma_{x_0}^1$ and therefore, by the continuity of the Neumann eigenvalues under $C^1$- perturbations, see \cite{Arrieta&Carvalho}, we have that $\lambda_2(\Gamma^1_{x_0})=0$. But this means that $\Gamma_{x_0}^1$ is not a connected domain and therefore $\Gamma_{x_0}^1$ is not diffeomorphic to the unit ball $B(0,1)$.

Hence, we obtain the first inequality with $\beta=\frac{1}{\hat{\lambda}_2}$.
%
%

For the inequality in the domain $Q_\varepsilon$, use the estimate in $Q$ and the appropriate change of variables in the integrals.  

%

\par\medskip\noindent ii)  Since $K\subset X_0^\alpha$ is a compact set, for $\eta>0$ there exist $u_0^1, ..., u_0^{k(\eta)}\in K$ such that 
$K\subset\bigcup_{i=1}^{k(\eta)} B(u_0^i, \eta). $
Then, for each $u_0\in K$, there exists $i\in\{1, 2, ..., k(\eta)\}$, such that $\|u_0^i-u_0\|_{X_0^\alpha}\leq\eta$.

Moreover, by the continuity of eigenvalues, see \cite[Section 3]{Arrieta-Santamaria-C0},  we have for each $i\in\{1,2, ..., k(\eta)\}$
\begin{equation}\label{convergence-i}
\|Eu_0^i\|_{X_\eps^\alpha}\rightarrow \|u_0^i\|_{X_0^\alpha}.
\end{equation}

If we write $Eu_0=E(u_0-u_0^i)+ Eu_0^i$. Then, from \eqref{Ealpha} we know that $\|E(u_0-u_0^i)\|_{X_\eps^\alpha}\leq \eta$.  This implies 
$$\left|\|Eu_0\|_{X_\eps^\alpha}-\|Eu_0^i\|_{X_\eps^\alpha}\right|\leq \eta.$$
From \eqref{convergence-i}, we know that there exists an $\eps(\eta)$ such that, for $0\leq\eps\leq\eps(\eta)$,
$$\left|\|Eu_0^i\|_{X_\eps^\alpha}-\|u_0^i\|_{X_0^\alpha}\right|\leq \eta, $$
and,
$$\left|\|Eu_0\|_{X_\eps^\alpha}\mathord-\|u_0^i\|_{X_0^\alpha}\right|=\left|\|\|Eu_0\|_{X_\eps^\alpha}\mathord-\|Eu_0^i\|_{X_\eps^\alpha}\mathord+\|Eu_0^i\|_{X_\eps^\alpha}\mathord-\|u_0^i\|_{X_0^\alpha}\right|\leq2\eta,$$
for $0\leq\eps\leq\eps(\eta)$.
So,
$$\left|\|Eu_0\|_{X_\eps^\alpha}-\|u_0\|_{X_0^\alpha}\right|=$$
$$\left|\|Eu_0\|_{X_\eps^\alpha}-\|Eu_0^i\|_{X_\eps^\alpha}\mathord+\|Eu_0^i\|_{X_\eps^\alpha}\mathord-\|u_0^i\|_{X_0^\alpha}\mathord+\|u_0^i\|_{X_0^\alpha}\mathord-\|u_0\|_{X_0^\alpha}\right|\leq 3\eta,$$
for $0<\eps\leq\eps(\eta)$. 

That is, for any $ K\subset X_0^\alpha$ and $K$ a compact set, 

$$\sup_{u_0\in K}\left|\|Eu_0\|_{X_\eps^\alpha}-\|u_0\|_{X_0^\alpha}\right|\rightarrow 0,\qquad\textrm{as}\quad\eps\rightarrow 0.$$
This concludes the proof. 
\begin{flushright}$\blacksquare$\end{flushright}

\section{Some previous results on convergence of Inertial Manifolds}
\label{previous}

In this section we are going to recall the results obtained in \cite{Arrieta-Santamaria-C0,Arrieta-Santamaria-C1} where we were able to analyze the convergence of inertial manifolds for abstract evolutionary equations under certain conditions. We were also able to obtain estimates on the distance of these inertial manifolds in the $C^0$ topology, see \cite{Arrieta-Santamaria-C0} and in the $C^{1,\theta}$ topology, see \cite{Arrieta-Santamaria-C1}.  We refer to these two papers for details.  

Hence, consider the family of abstract problems (like \eqref{problemaperturbado-thindomain}, \eqref{problemalimite-thindomain})
\begin{equation}\label{problemalimite}
(P_0)\left\{
\begin{array}{r l }
u^0_t+A_0u^0&=F_0^\eps(u^0),\\
u^0(0)\in X^\alpha_0,
\end{array}
\right.
\end{equation}
and 
\begin{equation}\label{problemaperturbado}
(P_\varepsilon)\left\{
\begin{array}{r l }
u^\varepsilon_t+A_\varepsilon u^\varepsilon&=F_\varepsilon(u^\varepsilon),\qquad 0<\varepsilon\leq \eps_0\\
u^\varepsilon(0)\in X^\alpha_\varepsilon,
\end{array}
\right.
\end{equation}
where we assume, that $A_\varepsilon$ is self-adjoint positive linear operator with compact resolvent on a separable real Hilbert space $X_\varepsilon$,  that is 
$A_\varepsilon: D(A_\varepsilon)=X^1_\varepsilon\subset X_\varepsilon\rightarrow X_\varepsilon,$
and  $F_\eps:X_\eps^\alpha\to X_\eps$, $F_0^\eps:X_0^\alpha\to X_0$ are nonlinearities guaranteeing  global existence of solutions of \eqref{problemaperturbado}, for each $0\leq \eps\leq \eps_0$ and for some $0\leq \alpha<1$.   Observe that for problem \eqref{problemalimite} we even assume that the nonlinearity depends on $\eps$ also. 


%

We also assume the existence of linear continuous operators, $E$ and $M$, such that, $E: X_0\rightarrow X_\varepsilon$, $M: X_\varepsilon\rightarrow X_0$ and $E_{\mid_{X^\alpha_0}}: X_0^\alpha\rightarrow X_\varepsilon^\alpha$ and $M_{\mid_{X_\varepsilon^\alpha}}: X_\varepsilon^\alpha\rightarrow X_0^\alpha$, satisfying,
\begin{equation}\label{cotaextensionproyeccion}
\|E\|_{\mathcal{L}(X_0, X_\varepsilon)}, \|M\|_{\mathcal{L}(X_\varepsilon, X_0)}\leq \kappa ,\qquad \|E\|_{\mathcal{L}(X^\alpha_0, X^\alpha_\varepsilon)}, \|M\|_{\mathcal{L}(X^\alpha_\varepsilon, X^\alpha_0)}\leq \kappa.
\end{equation}
for some constant $\kappa\geq 1$.   
We also assume these operators satisfy the following properties,
\begin{equation}\label{propiedadesextensionproyeccion}
M\circ E= I,\qquad \|Eu_0\|_{X_\varepsilon}\rightarrow \|u_0\|_{X_0}\quad\textrm{for}\quad u_0\in X_0.
\end{equation}

%

With respect to the relation between both operators, $A_0$ and $A_\eps$  and following \cite{Arrieta-Santamaria-C0,Arrieta-Santamaria-C1}, we will assume the following hypothesis

\vspace{0.5cm}
{\sl \paragraph{\textbf{(H1).}}  With $\alpha$ the exponent from problems (\ref{problemaperturbado}), we have
\begin{equation}\label{H1equation}
\|A_\varepsilon^{-1}- EA_0^{-1}M\|_{\mathcal{L}(X_\varepsilon, X_\varepsilon^\alpha)}\to 0\quad \hbox{ as } \eps\to 0.
\end{equation}
}
\par\bigskip 


Let us define $\tau(\eps)$ as an increasing function of $\eps$ such that 
\begin{equation}\label{definition-tau}
\|A_\varepsilon^{-1}E- EA_0^{-1}\|_{\mathcal{L}(X_0, X_\varepsilon^\alpha)}\leq \tau(\eps).
\end{equation}

\par\bigskip

We also recall hypothesis {\bf(H2)} from \cite{Arrieta-Santamaria-C0}, regarding the nonlinearities $F_0$ and $F_\eps$,   \par\bigskip
{\sl 
\paragraph{\textbf{(H2).}} We assume that the nonlinear terms  $F_\varepsilon: X^\alpha_\varepsilon\rightarrow X_\varepsilon$  and   $F_0^\varepsilon: X^\alpha_0\rightarrow X_0$ for $0< \eps\leq \eps_0$, satisfy:

\begin{enumerate}
\item[(a)]  They are uniformly bounded, that is, there exists a constant $C_F>0$ independent of $\varepsilon$ such that,
$$\|F_\varepsilon\|_{L^\infty(X_\varepsilon^\alpha, X_\varepsilon)}\leq C_F, \quad \|F_0^\varepsilon\|_{L^\infty(X_0^\alpha, X_0)}\leq C_F$$
\item[(b)] They are globally Lipschitz on $X^\alpha_\varepsilon$ with a uniform Lipstichz constant $L_F$, that is,
\begin{equation}\label{LipschitzFepsilon}
\|F_\varepsilon(u)- F_\varepsilon(v)\|_{X_\varepsilon}\leq L_F\|u-v\|_{X_\varepsilon^\alpha}
\end{equation}
\begin{equation}\label{LipschitzF0}
\|F_0^\varepsilon(u)- F_0^\varepsilon(v)\|_{X_0}\leq L_F\|u-v\|_{X_0^\alpha}. 
\end{equation}

\item[(c)] They have a uniformly bounded support for $0<\varepsilon\leq \eps_0$: there exists $R>0$ such that 
$$Supp F_\varepsilon\subset D_{R}=\{u_\varepsilon\in X_\varepsilon^\alpha: \|u_\varepsilon\|_{X_\varepsilon^\alpha}\leq R\}$$
$$Supp F_0^\varepsilon\subset D_{R}=\{u_0\in X_0^\alpha: \|u_0\|_{X_0^\alpha}\leq R\}.$$

\item[(d)]  $F_\eps$ is near $F_0^\eps$ in the following sense,
\begin{equation}\label{estimacionefes}
\sup_{u_0\in X^\alpha_0}\|F_\varepsilon (Eu_0)-EF_0^\eps (u_0)\|_{X_\varepsilon}=\rho(\varepsilon),
\end{equation}
and $\rho(\varepsilon)\rightarrow 0$ as $\varepsilon\rightarrow 0$.

\end{enumerate} 
}


%
%

The family of operators $A_\varepsilon$, for $0\leq \eps\leq \eps_0$,  are selfadjoint and have compact resolvent. This,  implies that their spectrum is discrete real and consists only of eigenvalues, each one with finite multiplicity. Moreover, the fact that $A_\varepsilon$, $0\leq \varepsilon\leq \eps_0$, is positive implies that its spectrum is positive. So, denoting by $\sigma(A_\varepsilon)$ the spectrum of the operator $A_\varepsilon$,  we have 
$$\sigma(A_\varepsilon)=\{\lambda_n^\varepsilon\}_{n=1}^\infty,\qquad\textrm{ and}\quad 0<c\leq\lambda_1^\varepsilon\leq\lambda_2^\varepsilon\leq...\leq\lambda_n^\varepsilon\leq...$$
We also denote by $\{\varphi_i^\varepsilon\}_{i=1}^\infty$ an associated orthonormal family of eigenfunctions,  by $\mathbf{P}_{\mathbf{m}}^{\bm\varepsilon}$ the canonical orthogonal projection onto the eigenfunctions, $\{\varphi^\varepsilon_i\}_{i=1}^m$, corresponding to the first $m$ eigenvalues of the operator $A_\varepsilon $, $0\leq\varepsilon\leq\varepsilon_0$ and by $\mathbf{Q}^{\bm\varepsilon}_{\mathbf{m}}$ the projetion over its orthogonal complement, see \cite{Arrieta-Santamaria-C0}. 
\par\bigskip

If we assume that {\bf (H1)} holds, then we obtain that the eigenvalues and eigenfunctions of the operator $A_\eps$ converge to the eigenvalues and eigenfunctions of $A_0$.  As a matter of fact, we get   that 
$$\lambda_i^\eps\eto \lambda_i^0, \hbox{ for each } i\in \N$$
and
\begin{equation}\label{convergence-of-projection}
\|\mathbf{P}_{\mathbf{m}}^{\bm\varepsilon}E(v)-E\mathbf{P}_{\mathbf{m}}^{0}(v)\|\leq C\tau(\eps)\|v\|_{X_0}
\end{equation}
(see \cite[Lemma 3.7]{Arrieta-Santamaria-C0}).

This last estimate implies easily that the set $$\{\mathbf{P}_{\mathbf{m}}^{\bm\varepsilon}(E\varphi^0_1), \mathbf{P}_{\mathbf{m}}^{\bm\varepsilon}(E\varphi^0_2), ...,\mathbf{P}_{\mathbf{m}}^{\bm\varepsilon}(E\varphi^0_m)\},\qquad\textrm{for}\quad 0\leq\varepsilon\leq\varepsilon_0,$$
 constitutes a basis in $\mathbf{P}_{\mathbf{m}}^{\bm\varepsilon}(X_\varepsilon)=[\psi_1^\varepsilon, ..., \psi_m^\varepsilon]$, that is, the space generated by the first $m$ eigenfunctions.

Let us denote by $j_\varepsilon$ the isomorphism from $ \mathbf{P}_{\mathbf{m}}^{\bm\varepsilon}(X_\varepsilon)=[\psi_1^\varepsilon, ..., \psi_m^\varepsilon]$ onto $\mathbb{R}^m$, that gives us the coordinates of each vector. That is,

\begin{equation}\label{definition-jeps}
\begin{array}{rl}
j_\varepsilon:\mathbf{P}_{\mathbf{m}}^{\bm\varepsilon}(X_\varepsilon)&\longrightarrow \mathbb{R}^m, \\
w_\varepsilon&\longmapsto\bar{p},
\end{array}
\end{equation}
where $w_\varepsilon=\sum^m_{i=1} p_i\psi^\varepsilon_i$ and $\bar{p}=(p_1, ..., p_m)$.

We denote by $|\cdot|$ the usual euclidean norm in $\mathbb{R}^m$, that is $|\bar{p}|=\left(\sum_{i=1}^mp_i^2\right)^{\frac{1}{2}}$, 
and by $|\cdot|_{\eps,\alpha}$ the following weighted one,
\begin{equation}\label{normaalpha}
|\bar{p}|_{\eps,\alpha}=\left(\sum_{i=1}^mp_i^2(\lambda_i^\varepsilon)^{2\alpha}\right)^{\frac{1}{2}}.
\end{equation}
\vspace{0.4cm}
We consider the spaces $(\mathbb{R}^m, |\cdot|)$ and $(\mathbb{R}^m, |\cdot|_{\eps,\alpha})$, that is, $\mathbb{R}^m$ with the norm $|\cdot|$ and $|\cdot|_{\eps,\alpha}$, respectively, and notice that for $w_0=\sum^m_{i=1} p_i\psi^0_i$ and $0\leq\alpha<1$ we have that,
\begin{equation}\label{normajepsilon}
\|w_0\|_{X^\alpha_0}=|j_0(w_0)|_{\eps,\alpha}.
\end{equation}

\par\bigskip 
We are looking for inertial manifolds for system \eqref{problemaperturbado} and \eqref{problemalimite} which will be obtained as graphs of appropriate functions.  This motivates the introduction of the sets
$\mathcal{F}_\eps(L,\rho)$ defined as 
$$\mathcal{F}_\eps(L,\rho)=\{ \Phi :\mathbb{R}^m\rightarrow\mathbf{Q}_{\mathbf{m}}^{\bm\varepsilon}(X^\alpha_\varepsilon),\quad\textrm{such that}\quad \textrm{supp } \Phi\subset B_R\quad \textrm{and}\quad$$
$$\quad \|\Phi(\bar{p}^1)-\Phi(\bar{p}^2)\|_{X^\alpha_\varepsilon}\leq L|\bar{p}^1-\bar{p}^2|_{\eps,\alpha} \quad\bar{p}^1,\bar{p}^2\in\mathbb{R}^m \}.$$ 
Then we can show the following result.
\begin{prop} (\cite{Arrieta-Santamaria-C0})\label{existenciavariedadinercial}
Let hypotheses {\bf (H1)} and {\bf (H2)} be satisfied. Assume also that $m\geq 1$ is such that,
\begin{equation}\label{CondicionAutovaloresFuerte0}
\lambda_{m+1}^0-\lambda_m^0\geq 3(\kappa+2)L_F\left[(\lambda_m^0)^\alpha+(\lambda_{m+1}^0)^\alpha\right],
\end{equation}
and
\begin{equation}\label{autovalorgrande0}
(\lambda_m^0)^{1-\alpha}\geq 6(\kappa +2)L_F(1-\alpha)^{-1}.
\end{equation}
Then, there exist $L<1$ and $\varepsilon_0>0$ such that for all $0<\varepsilon\leq\varepsilon_0$ there exist  inertial manifolds $\mathcal{M}_\varepsilon$ and $\mathcal{M}_0^\varepsilon$  for    (\ref{problemaperturbado}) and (\ref{problemalimite}) respectively,  given by the ``graph'' of a function $\Phi_\varepsilon\in\mathcal{F}_\eps(L,\rho)$ and $\Phi_0^\varepsilon\in\mathcal{F}_0(L,\rho)$. 
\end{prop}

\begin{re} 
%
We have written quotations in the word ``graph'' since the manifolds  $\mathcal{M}_\varepsilon$, $\mathcal{M}_0^\varepsilon$ are not properly speaking the graph of the functions $\Phi_\varepsilon$, $\Phi_0^\varepsilon$ but rather the graph of the appropriate function obtained via the isomorphism $j_\eps$ which identifies $\mathbf{P}_{\mathbf{m}}^{\bm\varepsilon}(X_\varepsilon^\alpha)$ with $\R^m$. 
That is, $\mathcal{M}_\eps=\{ j_\eps^{-1}(\bar p)+\Phi_\varepsilon(\bar p); \quad \bar p\in \R^m\}$ and 
$\mathcal{M}_0^\eps=\{ j_0^{-1}(\bar p)+\Phi_0^\varepsilon(\bar p); \quad \bar p\in \R^m\}$
\end{re}

The main result from \cite{Arrieta-Santamaria-C0} was the following:

\begin{teo} (\cite{Arrieta-Santamaria-C0})
\label{distaciavariedadesinerciales}
Let hypotheses {\bf (H1)} and {\bf (H2)} be satisfied and let $\tau(\eps)$ be defined by  \eqref{definition-tau}. 
Then, under the hypothesis of Proposition \ref{existenciavariedadinercial}, if $\Phi_\varepsilon$ are the maps  that give us the inertial manifolds for $0<\eps\leq \eps_0$ then we have,
\begin{equation}\label{distance-inertialmanifolds}
\|\Phi_\varepsilon-E\Phi_0^\eps\|_{L^\infty(\mathbb{R}^m, X^\alpha_\varepsilon)}\leq C[\tau(\varepsilon)|\log(\tau(\varepsilon))|+\rho(\varepsilon)],
\end{equation}
with $C$ a constant independent of $\varepsilon$.
\end{teo}
\par\bigskip

%
%

To obtain stronger convergence results on the inertial manifolds, we will need to requiere stronger conditions on the nonlinearites.  These conditions are stated in the following hypothesis, 
\par\bigskip
{\sl 
\paragraph{\textbf{(H2').}}
 We assume that the nonlinear terms $F_\varepsilon$ and $F_0^\eps$, satisfy hipothesis {\bf(H2)} and they are uniformly $C^{1,\theta_F}$ functions from $ X_\varepsilon^\alpha$ to $X_\varepsilon$, and $X_0^\alpha$ to $X_0$ respectively, for some $0<\theta_F\leq 1$. That is, $F_\eps\in C^1(X_\varepsilon^\alpha, X_\varepsilon)$, $F_0^\eps\in C^1(X_0^\alpha, X_0)$ and there exists  $L>0$, independent of $\eps$, such that 
$$\|DF_\varepsilon(u)-DF_\varepsilon(u')\|_{\mathcal{L}(X_\varepsilon^\alpha, X_\varepsilon)}\leq L\|u-u'\|^{\theta_F}_{X_\varepsilon^\alpha},\qquad \forall u, u'\in X_\varepsilon^\alpha.$$
$$\|DF_0^\varepsilon(u)-DF_0^\varepsilon(u')\|_{\mathcal{L}(X_0^\alpha, X_0)}\leq L\|u-u'\|^{\theta_F}_{X_0^\alpha},\qquad \forall u, u'\in X_0^\alpha.$$

}

We can show,
\begin{prop}\label{FixedPoint-E^1Theta} (\cite{Arrieta-Santamaria-C1})
 Assume hypotheses {\bf(H1)} and {\bf(H2')} are satisfied and that the gap conditions \eqref{CondicionAutovaloresFuerte0}, \eqref{autovalorgrande0} hold. Then, for any $\theta>0$ such that $\theta\leq\theta_F$ and $\theta<\theta_0$ (for certain $\theta_0>0$, see details in \cite{Arrieta-Santamaria-C1}) 
the functions $\Phi_\eps$, and $\Phi_0^\eps$ for  $0< \eps\leq\eps_0$, obtained above,  which give the inertial manifolds, are $C^{1, \theta}(\mathbb{R}^m, X_\eps^\alpha)$ and $C^{1, \theta}(\mathbb{R}^m, X_0^\alpha)$.  Moreover, the $C^{1,\theta}$ norm is bounded uniformly in $\eps$, for $\eps>0$ small. 
\end{prop}

The main result from \cite{Arrieta-Santamaria-C1} is the following:

\par\bigskip
\begin{teo}\label{convergence-C^1-theo} (\cite{Arrieta-Santamaria-C1})
Let hypotheses {\bf (H1)},  {\bf (H2')} and gap conditions \eqref{CondicionAutovaloresFuerte0}, 
\eqref{autovalorgrande0} be satisfied, so that we have inertial manifolds 
$\mathcal{M}^\eps$, $\mathcal{M}_0^\eps$ given as the graphs of the functions $\Phi_\eps$, $\Phi_0^\eps$ for $0<\eps\leq\eps_0$.   
If we denote by
\begin{equation}\label{convergenceDF}
\beta(\eps)=\sup_{u\in \mathcal{M}^\eps_0}\|DF_\varepsilon \big(Eu \big)E-EDF_0^\eps\big(u\big)\|_{\mathcal{L}(X_0^\alpha, X_\varepsilon)},
\end{equation}
%
%
%
%
%
%
%
 then, there exists $\theta^*$ with $0<\theta^*<\theta_F$ such that for all $0<\theta<\theta^*$, we obtain the following estimate
\begin{equation}\label{distance-C^1-inertialmanifolds}
\|\Phi_\varepsilon-E\Phi_0^\eps\|_{C^{1, \theta}(\mathbb{R}^m, X_\varepsilon^\alpha)}\leq \mathbf{C} \left(\left[\beta(\varepsilon)+\Big(\tau(\varepsilon)|\log(\tau(\varepsilon))|+\rho(\varepsilon)\Big)^{\theta^*}\right]\right)^{1-\frac{\theta}{\theta^*}},
 \end{equation}
where $\tau(\eps)$, $\rho(\varepsilon)$ are given by (\ref{definition-tau}), (\ref{estimacionefes}), respectively and $\mathbf{C}$ is a constant independent of $\eps$.
\end{teo}

\section{Estimates of the elliptic part}\label{elliptic}

As we mentioned in the introduction, a very important ingredient in comparing the dynamics of both problems is the convergence of the resolvent operators associated to the linear elliptic problems. In this section we will obtain these rates.

We consider the elliptic problems,

\begin{equation}\label{linearQ}
\left\{
\begin{array}{c l r}
-\frac{\partial^2u_\varepsilon}{\partial x^2}-\frac{1}{\varepsilon^2}\Delta_yu_\varepsilon+\mu u_\varepsilon\;&\; = f_\varepsilon, \;&\;\textrm{in}\quad Q\\
\frac{\partial u}{\partial\nu_x}+\frac{1}{\varepsilon^2}\frac{\partial u}{\partial\nu_\mathbf{y}}\;&\;=0\quad&\textrm{on}\quad \partial Q,
\end{array}
\right.
\end{equation} 
and 
\begin{equation}\label{linear01}
\left\{
\begin{array}{r l r}
-\frac{1}{g}(g {v_\varepsilon}_x)_x + \mu v_\varepsilon\;&\; = h_\eps, \;&\;\textrm{in}\quad (0, 1)\\
v_{\varepsilon_x}(0)=v_{\varepsilon_x}(1)\;&\; =0,\;&\;
\end{array}
\right.
\end{equation} 
with $f_\varepsilon\in L^2(Q)$, $u_\varepsilon\in H^1_{\bm\varepsilon}(Q)$ and $h_\eps\in L^2_g(0,1)$, $v_\varepsilon\in H_g^1(0,1)$.  Notice that the existence and uniqueness of solutions of the problems above is guaranteed by Lax-Milgram theorem. 

\par\bigskip 

We can prove the following key result.
\begin{prop}\label{resolvente}
Let  $f_\varepsilon\in L^2(Q)$ and let $h_\eps=M f_\eps$. We define the functions $u_\varepsilon\in H^1_{\bm\varepsilon}(Q)$ and $v_\varepsilon\in H^1_g(0,1)$ as the solutions of the linear problems (\ref{linearQ}) and (\ref{linear01}), respectively. Then, there exist a constant $C>0$ independent of $\eps$ and $f_\varepsilon$ such that, 
$$\|u_\varepsilon-Ev_\varepsilon\|_{H^1_{\bm\varepsilon}(Q)}\leq C\varepsilon\|f_\varepsilon\|_{L^2(Q)}.$$
\end{prop}

\paragraph{\sl Proof.}  Since the proof of this result is technical, we prefer to postpone its proof for later. We provide the proof of this result in Appendix \ref{proof-resolvente}.

\vspace{0.8cm}

\begin{re}\label{estimaciondebilresolvente}
Note that if we consider problems \eqref{linearQ} and \eqref{linear01} with $f_\varepsilon=E h_\varepsilon$ then, $\|f_\varepsilon-M_\varepsilon f_\varepsilon\|_{L^2(Q_\varepsilon)}=0$ and so, we obtain the same estimate,
\begin{equation}\label{estimaciondebilresolvente-0}
\|u_\varepsilon-Ev_\varepsilon\|_{H^1_{\bm\varepsilon}(Q)}\leq C\varepsilon\|Eh_\varepsilon\|_{L^2(Q)}.
\end{equation}
\end{re}

\begin{re}\label{version-abstracta}
Writing this proposition in the abstract setting we get 
\begin{equation}\label{resolventefuerte}
\|A_\varepsilon^{-1}-EA_0^{-1}M\|_{\mathcal{L}(L^2(Q), H^1_{\bm\varepsilon}(Q))}\leq C\varepsilon,
\end{equation}
Observe that \eqref{resolventefuerte} implies that hypothesis {\bf (H1)} holds for $\alpha=1/2$ and therefore for $0\leq \alpha\leq 1/2$. Moreover, estimate \eqref{estimaciondebilresolvente-0} shows that $\tau(\eps)$ satisfies $\tau(\eps)=\eps$.  
\end{re}

%

\par\bigskip 
We show now, in a formal way, that the estimate obtained in Proposition \ref{resolvente} is optimal. For this, we will consider a domain $Q_\eps$ having circular cross sections and with the aid of an asymptotic expansion of the solution $u_\eps$, we will obtain that the estimates obtained are optimal. 

Hence, let $Q_\varepsilon=\{(x, \varepsilon\mathbf{y})\in\mathbb{R}^d : (x, \mathbf{y})\in Q\}$, with $\varepsilon\in (0, 1)$, and $Q=\{(x, \mathbf{y})\in\mathbb{R}^d : 0\leq x\leq 1, |\mathbf{y}|< r(x) \}$, so that the transversal sections $\Gamma_x^1$ of the domain $Q$ are disks  centered at the origin of radius $r(x)$. Obviously, the change of variables which takes $Q_\varepsilon$ into the fixed domain $Q$ is the following,
$$X=x, \qquad \mathbf{Y}=\varepsilon \mathbf{y},$$
with $(X, \mathbf{Y})\in Q_\eps$ and $(x, \mathbf{y})\in Q$.
This change of variables transforms the original problem into the following linear problem in $Q$ (we consider the coefficient of equation $\alpha=1$),  
\begin{equation}\label{linearQexample}
\left\{
\begin{array}{r l r}
-\frac{\partial^2u_\varepsilon}{\partial x^2}-\frac{1}{\varepsilon^2}\Delta_{\mathbf{y}}u_\varepsilon+u_\varepsilon\;&\; = Ef, \;&\;\textrm{in}\quad Q\\
(\nabla_xu_\varepsilon, \frac{1}{\varepsilon}\nabla_{\mathbf{y}}u_\varepsilon)\cdot\nu_\varepsilon\;&\;=0\quad&\textrm{on}\quad \partial Q,
\end{array}
\right.
\end{equation} 
with $\partial Q=\Gamma_0^1\cup\partial_lQ\cup\Gamma_1^1$, where $\partial_l Q$ is the ``lateral boundary" which is given by 
$\partial_l Q=\{ (x,\mathbf{y}):   \mathbf{y}\in \partial \Gamma_x^1\}$ and 
\begin{equation}
\nu_\varepsilon=
\left\{
\begin{array}{c r}
(-1, 0),\qquad & \textrm{in}\quad \Gamma_0^1,\\
\left(\frac{-\varepsilon r r'}{r\sqrt{\varepsilon^2{r'}^2+1}}, \frac{\mathbf{y}}{r\sqrt{\varepsilon^2{r'}^2+1}}\right),\qquad & \textrm{in}\quad \partial \Gamma_x^1,\\
(1, 0),\qquad & \textrm{in}\quad \Gamma_1^1.
\end{array}
\right.
\end{equation}
The limit problem is given by 
\begin{equation}\label{linear01example}
\left\{
\begin{array}{r l r}
-\frac{1}{g}(g {v_0}_x)_x + v_0\;&\; = f, \;&\;\textrm{in}\quad (0, 1)\\
{v_0}_x(0)\;&\;=0,\;&\;{v_0}_x(1)=0,
\end{array}
\right.
\end{equation} 
with $f\in L^2_g(0,1)$.  Recall that  $g(x)=|\Gamma_x^1|=r(x)^{d-1}\omega_{d-1}$ and $\omega_{d-1}$ is the $(d-1)$-measure of the unit ball in $\mathbb{R}^{d-1}$.   

To analyze the rate of convergence of $u_\varepsilon\rightarrow v_0$ as $\varepsilon\rightarrow 0$, we express the solution of (\ref{linearQexample}) as the series
$$u_\varepsilon=\sum_{i=0}^\infty \varepsilon^iV_i(x, \mathbf{y})=V_0(x,\mathbf{y})+\eps V_1(x,\mathbf{y})+\eps^2 V_2(x,\mathbf{y})+\ldots$$
Introducing this expression in problem (\ref{linearQexample}) we obtain,
\begin{equation}\label{problemserie}
\left\{
\begin{array}{r l r}
-\sum_{i=0}^\infty \varepsilon^iV_{i_{xx}}-\frac{1}{\varepsilon^2}\sum_{i=0}^\infty \varepsilon^i \Delta_\mathbf{y} V_{i}+\sum_{i=0}^\infty \varepsilon^iV_i(x, \mathbf{y})\;&\; = Ef, \;&\;\textrm{in}\quad Q\\ \\
-\sum_{i=0}^\infty \varepsilon^iV_{i_x}\;&\;= 0,  \;&\;\textrm{on}\quad \Gamma_0^1\\ \\
-\sum_{i=0}^\infty \varepsilon^{i+1}V_{i_x}rr'+\sum_{i=0}^\infty \varepsilon^{i-1}\nabla_{\mathbf{y}} V_{i}\cdot \mathbf{y}\;&\;=0\quad&\textrm{on}\quad \partial \Gamma_x^1\\ \\
\sum_{i=0}^\infty \varepsilon^iV_{i_x}\;&\;= 0,\quad&\textrm{on}\quad  \Gamma_1^1 
\end{array}
\right.
\end{equation} 
\par\bigskip
Putting in groups of powers of $\varepsilon$, we have the following equalities in $Q$,
\begin{equation}\label{Equation-in-Q}
\begin{array}{c}
\Delta_{\mathbf{y}}V_0(x, \mathbf{y})=0, \\ \\
\Delta_{\mathbf{y}}V_1(x, \mathbf{y})=0, \\ \\
-V_{0_{xx}}(x, \mathbf{y})-\Delta_{\mathbf{y}}V_{2}(x, \mathbf{y})+V_0(x, \mathbf{y})-f(x)=0, \\ \\
-V_{i_{xx}}(x, \mathbf{y})-\Delta_{\mathbf{y}}V_{i+2}(x, \mathbf{y})+V_i(x, \mathbf{y})=0,\qquad \textrm{for}\quad i=1, 2,...
\end{array}
\end{equation}

and, from the boundary condition, we have,
\begin{equation}\label{Equation-in-Boundary-Q}
\begin{array}{c}
V_{i_x}(x, \mathbf{y})=0\qquad\textrm{on}\quad \Gamma_0^1\cup\Gamma_1^1,\qquad\textrm{for}\quad i=0, 1, 2, ...\\ \\
\nabla_{\mathbf{y}}V_0(x, \mathbf{y})\cdot\mathbf{y}=0,  \quad \hbox{ on  } \partial \Gamma_x^1 \\ \\
\nabla_{\mathbf{y}}V_1(x, \mathbf{y})\cdot\mathbf{y}=0, \quad \hbox{ on  } \partial \Gamma_x^1\\ \\
-V_{i_x}(x, \mathbf{y})rr'+\nabla_{\mathbf{y}}V_{i+2}(x, \mathbf{y})\cdot\mathbf{y}=0,\qquad\textrm{for}\quad i=0, 1, 2,...\quad \hbox{ on  } \partial \Gamma_x^1
\end{array}
\end{equation}

First, for $x\in(0, 1)$ fixed, we focus in the particular problems in $\Gamma_x^1$ in which $V_0(x, \mathbf{y})$ and $V_1(x, \mathbf{y})$ are involved,
\begin{equation}\label{problemV0V1}
\begin{array}{r l r c r l r}
\Delta_{\mathbf{y}}V_0(x,\mathbf{y})\;&\;=0\quad& \textrm{in}\quad \Gamma_x^1, &\hspace{1cm}& \Delta_{\mathbf{y}}V_1(x,\mathbf{y})\;&\;=0\qquad& \textrm{in}\quad \Gamma_x^1\\
\nabla_{\mathbf{y}}V_0(x, \mathbf{y})\cdot\mathbf{y}\;&\;=0\quad& \textrm{on}\quad \partial\Gamma_x^1, &\hspace{1cm}&\nabla_{\mathbf{y}}V_1(x, \mathbf{y})\cdot\mathbf{y}\;&\;=0\qquad& \textrm{on}\quad \partial\Gamma_x^1.
\end{array}
\end{equation} 
Both problems imply that, for each $x\in(0, 1)$, $V_0(x, \mathbf{y})$ and $V_1(x, \mathbf{y})$ are constant in $\Gamma_x^1$. It means both functions only depend on $x$,
$$V_0(x, \mathbf{y})= V_0(x), \hspace{2cm} V_1(x, \mathbf{y})=V_1(x).$$
Since $V_0$ only depends on $x$, the third condition in \eqref{Equation-in-Q}  and in \eqref{Equation-in-Boundary-Q} can
be written as 
\begin{equation}\label{problemV2}
\left\{
\begin{array}{r l r}
\Delta_{\mathbf{y}}V_2(x,\mathbf{y})\;&\;=-V_{0_{xx}}(x)+V_0(x)-f(x)\qquad& \textrm{in}\quad \Gamma_x^1, \\
\nabla_{\mathbf{y}}V_2(x, \mathbf{y})\cdot\nu\;&\;=V_{0_x}(x)r'(x)\qquad& \textrm{on}\quad \partial\Gamma_x^1.
\end{array}
\right.
\end{equation} 
Integrating over $\Gamma_1^x$ in the equation and using the boundary condition, we find that in order to have solutions of \eqref{problemV2} we must have (Fredholm alternative), 
$${V_0}_xr'|\partial\Gamma_x^1|=  (-V_{0_{xx}}(x)+V_0(x)-f(x))|\Gamma_x^1|.$$
That is,
$$-V_{0_{xx}}(x)+V_0(x)-f(x)={V_0}_xr'  \frac{|\partial\Gamma_x^1|}{|\Gamma_x^1|}=\frac{d-1}{r}V_{0_x}r'.$$
%
%
Now, since $\frac{g_x}{g}=(d-1)\frac{r'}{r}$ we easily get
$$-\frac{1}{g}(gV_{0_{x}})_x+V_0=f(x),$$
and the boundary conditions are given by 
$$V_{0_x}(0)=V_{0_x}(1)=0.$$

This implies $V_0(x, \mathbf{y})=v_0(x)$ is the solution of the limit problem \eqref{linear01example}.  Moreover, the function $V_2(x,\mathbf{y})$ satisfies \eqref{problemV2} and it is not identically 0 in general (if for instance $f\not\equiv 0$). 

\par\bigskip Proceeding in a similar way with $V_1$ and $V_3$ we get, 
\begin{equation}\label{problemV3constant}
\left\{
\begin{array}{r l r}
\Delta_{\mathbf{y}}V_3(x,\mathbf{y})\;&\;=-V_{1_{xx}}(x)+V_1(x)\qquad& \textrm{in}\quad \Gamma_x^1, \\
\nabla_{\mathbf{y}}V_3(x, \mathbf{y})\cdot\nu\;&\;=r'V_{1_x}(x)\qquad& \textrm{on}\quad \partial\Gamma_x^1,
\end{array}
\right.
\end{equation} 
and with the Fredholm alternative, the function $V_1$ needs to satisfy $-\frac{1}{g}(gV_{1_{x}})_x+V_1=0,$ with the boundary conditions ${V_1}_x(0)={V_1}_x(1)=0$  (see \eqref{Equation-in-Boundary-Q}).   This implies that $V_1(\cdot)\equiv 0$ and from 
\eqref{problemV3constant} we get $V_3=V_3(x)$.  With an induction argument it is not difficult to see now that $V_i\equiv 0$ for all odd $i$.   Hence, 
$$u_\eps(x,\mathbf{y})=v_0(x)+\eps^2V_2(x,\mathbf{y})+\eps^4 V_4(x,\mathbf{y}) +\ldots$$
where $V_2(x,\mathbf{y})$ is the solution of \eqref{problemV2} which is generically non zero. 
%
%
%

Then, for $\varepsilon$ small enough,
$$\|u_\varepsilon-EV_0\|_{H^1_{\bm\varepsilon}(Q)}=\|\varepsilon^2V_2(x, \mathbf{y})\|_{H^1_{\bm\varepsilon}(Q)}+\ldots$$
$$=\left(\varepsilon^4\int_Q(|\nabla_xV_2|^2+\frac{1}{\varepsilon^2}|\nabla_{\mathbf{y}}V_2|^2+|V_2|^2)dxd{\mathbf{y}}\right)^{\frac{1}{2}}+\ldots = \eps\|\nabla_\mathbf{y}V_2\|_{L^2(Q)}+o(\eps).$$
But, 
$$\|\nabla_\mathbf{y}V_2\|_{L^2(Q)}\sim \|f\|_{L^2(Q)},$$
which implies that estimate from Proposition \ref{resolvente} is optimal. 
%
%

\section{Analysis of the nonlinear terms}\label{nonlinear}
In this section we focus our study in the nonlinear terms. We will analyze its differentiability properties and we will prepare the nonlinearities to apply the results on existence and convergence of inertial manifolds described in Section \ref{previous} (see also \cite{Arrieta-Santamaria-C0,Arrieta-Santamaria-C1}).  As  a matter of fact, we will show that with appropriate cut-off functions the new nonlinearities satisfy hypotheses {\bf (H2)} and {\bf (H2')}  easing our way to the construction of the inertial manifolds and to estimating the distance between them.

\medskip

Recall that we denote by  $X_\eps^\alpha$, for $0\leq\eps\leq\eps_0$ the fractional power spaces corresponding to the elliptic operators, see Section \ref{setting}.

\medskip

First, we analyze the properties the nonlinear terms satisfy. Remember that the nonlinearity $f$, together with its first and second derivative  satisfy the boundedness condition (\ref{dissipativecondition2}). We denote by $F_\varepsilon: X_\eps^\alpha\rightarrow L^2(Q)$ the Nemytskii operator corresponding to $f$, that is,
\begin{equation}\label{Nem-op}
\begin{array}{rl}
F_\varepsilon:X_\eps^\alpha&\longrightarrow L^2(Q), \\
u&\longmapsto f(u),
\end{array}
\end{equation}
 \begin{equation}
\begin{array}{rl}
F_0:X_0^\alpha&\longrightarrow L^2_g(0, 1),\\
u&\longmapsto f(u).
\end{array}
\end{equation}
\bigskip

Then we have the following result.
\begin{lem}\label{nonlinearity-C1}
The Nemytskii operator $F_\varepsilon$, $\eps\geq 0$, satisfies the following properties:
\begin{itemize}
\item[(i)]$F_\eps$ is uniformly bounded from $X_\eps^\alpha$ into $L^2(Q)$. That is, there exists a constant $C_F>0$ independent of $\eps$ such that,
$$\|F_\eps\|_{L^\infty(X_\eps^\alpha, L^2(Q))}\leq C_F.$$
\item[(ii)]There exists $\theta_F\in (0,1]$ such that $F_\eps$ is $C^{1, \theta_F}(X_\eps^\alpha, L^2(Q))$ uniformly in $\eps$. That is, there exists a constant $L_F>0$, such that,
$$\|F_\eps(u)-F_\eps(v)\|_{L^2(Q)} \leq L_F\|u-v\|_{X_\eps^\alpha},$$
$$\|DF_\eps(u)-DF_\eps(v)\|_{\mathcal{L}(X_\eps^\alpha, L^2(Q))} \leq L_F\|u-v\|^{\theta_F}_{X_\eps^\alpha}$$
for all $u,v\in X_\eps^\alpha$ and all $0\leq \eps\leq \eps_0$: 
\end{itemize}
\end{lem}
\par\bigskip
\paragraph{\sl Proof. }
Item (i) is directly proved as follows. Since nonlinearity $f$ is uniformly bounded, see (\ref{dissipativecondition2}), 
$$\|F_\eps\|_{L^\infty(X_\eps^\alpha, L^2(Q))}=\sup_{u\in X_\eps^\alpha}\left(\int_Q|f(u(x, \mathbf{y}))|^2dxd\mathbf{y}\right)^{\frac{1}{2}}\leq L_f|Q|^{\frac{1}{2}},$$
for any $\eps\geq 0$ and $|Q|$ the Lebesgue measure of $Q$. So, we have the desired estimate with $C_F=L_f|Q|^{\frac{1}{2}}$.

To prove item (ii), we proceed as follows.
$$\|F_\eps(u)- F_\eps(v)\|_{L^2(Q)}=\left(\int_Q|f(u(x, \mathbf{y}))-f(v(x, \mathbf{y}))|^2dxd\mathbf{y}\right)^{\frac{1}{2}}.$$
Since $f$ is globally Lipschitz, see \eqref{dissipativecondition2},  then,
$$\|F_\eps(u)- F_\eps(v)\|_{L^2(Q)}\leq L_f\left(\int_Q|u(x, \mathbf{y})-v(x, \mathbf{y})|^2dxd\mathbf{y}\right)^{\frac{1}{2}}=$$
$$=L_f\|u-v\|_{L^2(Q)}\leq L_f\|u-v\|_{X_\eps^\alpha},$$
taking $L_F=L_f$ we have, for all $\eps\geq 0$, that $F_\eps$ is globally Lipschitz from $X_\eps^\alpha$ into $L^2(Q)$ with uniform constant $L_F$. To show the remaining part, notice first that for $u\in X_\eps^\alpha$, $DF_\eps(u)$ is given by the operator
\begin{equation}
\begin{array}{rl}
DF_\eps(u): X_\eps^\alpha&\longrightarrow L^2(Q), \\
v&\longmapsto f'(u)v,
\end{array}
\end{equation}
which is easily shown from the definition of Fr\'{e}chet derivative, the Sobolev embeddings $X_\eps^\alpha\hookrightarrow L^q$ for $q>2$, and the property \eqref{dissipativecondition2}. That is,
$$\|F_\eps(u+v)-F_\eps(u)-f'(u)v\|_{L^2(Q)}=\|\left(f'(\xi)-f'(u)\right)v\|_{L^2(Q)},$$
with $\xi$ an intermediate point between $u$ and $u+v$.

But, by \eqref{dissipativecondition2} $|f'(\xi)-f'(u)|\leq 2 L_f$ and also by the mean value theorem $|f'(\xi)-f'(u)|\leq L_f|\xi-u|\leq L_f |v|$. This implies  $\left|f'(\xi)-f'(u)\right|\leq 2 L_f |v|^\theta$, for all $0<\theta<1$.

Hence,
$$\|F_\eps(u+v)-F_\eps(u)-f'(u)v\|_{L^2(Q)}\leq 2L_f\|v^{1+\theta}\|_{L^2(Q)}=2L_f \|v\|_{L^{2+2\theta}(Q)}^{1+\theta}.$$
Choosing $2+2\theta < q$ we get that $DF_\eps(u)v=f'(u)v$.
\bigskip

 Moreover, we have that, for all $\eps\geq 0$,
$$\|DF_\eps(u)-DF_\eps(v)\|_{\mathcal{L}(X_\eps^\alpha, L^2(Q))}=\sup_{\phi\in X_\eps^\alpha,\,\, \|\phi\|_{X_\eps^\alpha}\leq 1}\|DF_\eps(u)\phi-DF_\eps(v)\phi\|_{L^2(Q)}.$$
Hence,
$$\|DF_\eps(u)-DF_\eps(v)\|_{\mathcal{L}(X_\eps^\alpha, L^2(Q))}=\sup_{\phi\in X_\eps^\alpha,\,\, \|\phi\|_{X_\eps^\alpha}\leq 1}\left(\int_Q(f'(u)-f'(v))^2\phi^2dxd\mathbf{y}\right)^{\frac{1}{2}}.$$
Note that, by H\"{o}lder inequality with exponents $\frac{d}{4\alpha}$ and $\frac{d}{d-4\alpha}$, (remember $\alpha<\frac{1}{2}$ and $d\geq 2$, so that both $\frac{d}{4\alpha}$, $\frac{d}{d-4\alpha} \in (1,\infty)$), we have,
$$\int_Q(f'(u)-f'(v))^2\phi^2dxd\mathbf{y}\leq \left(\int_Q|f'(u)-f'(v)|^{\frac{d}{2\alpha}} dxd\mathbf{y}\right)^{\frac{4\alpha}{d}}\left(\int_Q |\phi|^{\frac{2d}{d-4\alpha}}dxd\mathbf{y}\right)^{\frac{d-4\alpha}{d}}.$$

Then, from Remark \ref{EMalpha} ii)   we have,
$$\int_Q(f'(u)-f'(v))^2\phi^2dxd\mathbf{y}\leq C\left(\int_Q|f'(u)-f'(v)|^{\frac{d}{2\alpha}} dxd\mathbf{y}\right)^{\frac{4\alpha}{d}}\|\phi\|^2_{X_\eps^\alpha}. $$
Then,
$$\sup_{\phi\in X_\eps^\alpha,\,\, \|\phi\|_{X_\eps^\alpha}\leq 1}\left(\int_Q(f'(u)-f'(v))^2\phi^2dxd\mathbf{y}\right)^{\frac{1}{2}}\leq\left(\int_Q|f'(u)-f'(v)|^{\frac{d}{2\alpha}} dxd\mathbf{y}\right)^{\frac{2\alpha}{d}}$$
Next, note that, on the one side, by the mean value theorem and using \ref{dissipativecondition2}, we have,
$$|f'(u)-f'(v)|\leq L_f|u-v|.$$
On the other side, again by (\ref{dissipativecondition2}),
$$|f'(u)-f'(v)|\leq 2L_f.$$
Hence,
$$|f'(u)-f'(v)|\leq 2L_f\min\{1, |u-v|\}\leq 2L_f|u-v|^\theta,$$
for any $0\leq\theta\leq 1$, where we have used that if $0\leq x\leq 1$  and $0\leq \theta\leq 1$ then $x\leq x^\theta$. 
Then,
$$\|DF_\eps(u)-DF_\eps(v)\|_{\mathcal{L}(X_\eps^\alpha, L^2(Q))}\leq 2L_f\left(\int_Q|u-v|^{\frac{\theta d}{2\alpha}}dxd\mathbf{y}\right)^{\frac{2\alpha}{d}}.$$
Taking $\theta_F=\min\{1,\frac{4\alpha}{d-4\alpha}\}$,
$$\|DF_\eps(u)-DF_\eps(v)\|_{\mathcal{L}(X_\eps^\alpha, L^2(Q))}\leq 2L_f\left(\int_Q|u-v|^{\frac{2d}{d-4\alpha}}dxd\mathbf{y}\right)^{\frac{2\alpha}{d}}=2L_f\|u-v\|^{\theta_F}_{L^{\frac{2d}{d-4\alpha}}(Q)}.$$
Applying again the uniform embedding described in Remark \ref{EMalpha} ii), we obtain
$$\|DF_\eps(u)-DF_\eps(v)\|_{\mathcal{L}(X_\eps^\alpha, L^2(Q))}\leq 2L_f\|u-v\|_{X_\eps^\alpha}^{\theta_F}. $$
Taking $L_F=2L_f$ we have the result.
\begin{flushright}$\blacksquare$\end{flushright}

\begin{re}\label{alpha-positive}
i) Note that we have to impose $\alpha$ strictly positive to guarantee the smoothness of $F_\eps$, that is, to ensure that $F_\eps\in C^{1, \theta}(X_\eps^\alpha, L^2(Q))$ for $\theta$ small enough. As a matter of fact if $\alpha=0$ $X_\eps^\alpha=L^2(Q)$, any nonlinearity $F:L^2(Q)\to L^2(Q)$ which is a Nemytskii operator, as in \eqref{Nem-op}, cannot be $C^1$, unless it is linear, see \cite{Henry1}, Exercise 1. Although in \cite{Henry1}, the author considers the case $F_\eps(u)=sin (u)$, the argument can be easily extended to any $C^2$ function. 
\par ii) if $d\geq 4$ we always have that $\theta_F=\frac{4\alpha}{d-4\alpha}<1$ because $\alpha<1/2$.  Only in dimensions $d=2,3$ and choosing $\alpha<1/2$ but close enough to $1/2$ we may get $\frac{4\alpha}{d-4\alpha}>1$ and therefore $\theta_F=1$.  As a matter of fact, in dimensions $d=2,3$ we may show some higher differentiability of $F$. 
\end{re}
\bigskip

We fix $\alpha$ with $0<\alpha<\frac{1}{2}$.

\bigskip

As we have mentioned above, one of our basic tools consists in constructing inertial manifolds to reduce our problem to a finite dimensional one. In order to construct these manifolds and following \cite{Sell&You}, we need to ``prepare" the non-linear term making an appropriate cut off of the nonlinearity in the $X_\eps^\alpha$ norm, as it is done in \cite{Sell&You} . 

Next, we proceed to introduce this cut off.  For this, we start considering a function $\hat\Theta: \mathbb{R}\to [0,1]$ which is $C^\infty$ with compact support and such that  

\begin{equation}\label{funcionthetaR}
\hat\Theta(x)=
\left\{
\begin{array}{l c r}
1 & if & |x|\leq R^2\\
0 & if & |x|\geq4R^2.
\end{array}
\right.
\end{equation}
for some $R>0$, which in general will be large enough.  We will denote this function $\hat\Theta^R(x)$ if we need to make explicit its dependence on the parameter $R$.   With this function we define now
$\Theta_\eps:X_\eps^\alpha\to \mathbb{R}$ as $\Theta_\eps(u)=\hat \Theta(\|u\|^2_{X_\eps^\alpha})$ for $0\leq\eps\leq\eps_0$, and observe that $\Theta_\eps(u)=1$ if $\|u\|_{X_\eps^\alpha}\leq R$ and $\Theta_\eps(u)=0$ if $\|u\|_{X_\eps^\alpha}\geq 2R$ and again we will denote $\Theta_\eps$ by $\Theta_\eps^R$ if we need to make explicit its dependence on $R$. 

Now,  for $R>0$, large enough, and $0<\varepsilon\leq\eps_0$, we introduce the new nonlinear terms
\begin{equation}\label{definition-cutoff-eps}
\tilde{F}_\varepsilon(u_\eps):= \Theta_\varepsilon^R(u_\eps) F_\varepsilon(u_\eps),
\end{equation}
\begin{equation}\label{definition-cutoff-*}
\tilde{F}^\eps_0(u_0):= \Theta_\varepsilon^R(Eu_0) F_0(u_0),
\end{equation}
and
\begin{equation}\label{definition-cutoff-0}
\tilde{F}_0(u_0):= \Theta_0^R(u_0) F_0(u_0),
\end{equation}
We replace $F_\varepsilon$ and $F_0$ with the new nonlinearities $\tilde{F}_\varepsilon$, $\tilde{F}^\eps_0$ and $\tilde{F}_0$. Hence, now we have three systems, two of them in the limit space $X_0^\alpha$,
\begin{equation}\label{system-eps}
{u}_t=-A_\eps u+\tilde{F}_\eps(u),\qquad u\in X_\eps^\alpha
\end{equation}

\begin{equation}\label{system-*}
{u}_t=-A_0u+\tilde{F}^\eps_0(u),\qquad u\in X_0^\alpha,
\end{equation}

\begin{equation}\label{system-0}
{u}_t=-A_0u+\tilde{F}_0(u), \qquad u\in X_0^\alpha.
\end{equation}

Note that, since systems \eqref{system-*} and \eqref{system-0} share the linear part and $\tilde{F}_0(u)=\tilde{F}^\eps_0(u)$ for $\|u\|_{X_0^\alpha}\leq R$, then the attractor related to \eqref{system-*} and \eqref{system-0} coincides and it is $\mathcal{A}_0$. Moreover, although $\,\tilde{F}^\eps_0, \tilde{F}_0: X_0^\alpha\rightarrow X_0\,$, the nonlinearity $\tilde{F}^\eps_0$ depends on $\eps$.

\bigskip

\begin{re}
It may sound somehow strange the need to consider now three systems instead of the natural two (the perturbed one \eqref{system-eps} and the completely unperturbed one \eqref{system-0}).   The three systems meet the conditions to have inertial manifolds and we will see that they all are nearby in the $C^{1}$ topology. But, as we will see below, we 
will have good estimates for the distance between  the inertial manifolds for systems \eqref{system-eps} and \eqref{system-*} but not so good estimates for the distance between the inertial manifolds for systems  \eqref{system-eps} and \eqref{system-0} or 
 \eqref{system-*} and \eqref{system-0}.   
%
%
\end{re}

\bigskip

First, we analyze the properties $\tilde{F}_\eps$, $\tilde{F}^\eps_0$ and $\tilde{F}_0$ satisfy.
\begin{lem}\label{nonlinearity}
Let $\tilde{F}_\eps$, $0<\eps\leq\eps_0$, $\tilde{F}^\eps_0$ and $\tilde{F}_0$, be the new nonlinearities described above. Then they satisfy the following properties:

\begin{itemize}
\item[(a)]$\tilde{F}_\varepsilon(u)=F_\varepsilon(u)$, for all $u\in X_\eps^\alpha$, such that $\|u\|_{ X_\eps^\alpha}\leq R$, $\varepsilon> 0$ and $\tilde{F}^\eps_0(u_0)=F_0(u_0)$, $\tilde{F}_0(u_0)=F_0(u_0)$, for all  $u_0\in X_0^\alpha$, such that $\|Eu_0\|_{X_\eps^\alpha}\leq R$ and $\|u_0\|_{X_0^\alpha}\leq R$, respectively.
\item[(b)]$\tilde{F}_\varepsilon$ is $C^{1, \theta_F}(X_\eps^\alpha, L^2(Q))$ and $\tilde{F}^\eps_0$, $\tilde{F}_0$ are $C^{1, \theta_F}(X_0^\alpha, L_g^2(0,1))$  with $\theta_F$ the one from Lemma \ref{nonlinearity-C1}. That is, they are globally Lipschitz from $X_\eps^\alpha$ to $L^2(Q)$ and from $X_0^\alpha$ to $L_g^2(0,1)$, we denote by $L_F$ their Lipschitz constant, and 
\begin{equation}\label{Holder}
\|D\tilde{F}_\eps(u)-D\tilde{F}_\eps(u')\|_{\mathcal{L}(X_\eps^\alpha, L^2(Q))}\leq L_F\|u-u'\|_{X_\eps^\alpha}^{\theta_F},
\end{equation}

\begin{equation}\label{Holder*}
\|D\tilde{F}^\eps_0(u)-D\tilde{F}^\eps_0(u')\|_{\mathcal{L}(X_0^\alpha, L_g^2(0,1))}\leq L_F\|u-u'\|_{X_0^\alpha}^{\theta_F},
\end{equation}
\begin{equation}\label{Holder0}
\|D\tilde{F}_0(u)-D\tilde{F}_0(u')\|_{\mathcal{L}(X_0^\alpha, L_g^2(0,1))}\leq L_F\|u-u'\|_{X_0^\alpha}^{\theta_F},
\end{equation}
with $L_F$ independent of $\eps$.
\item[(c)]They are uniformly bounded,
$$\|\tilde{F}_\eps\|_{L^\infty(X_\eps^\alpha, L^2(Q))}\leq C_F,\,\,\,\,\,\,\,\,\,\|\tilde{F}^\eps_0\|_{L^\infty(X_0^\alpha, L_g^2(0,1))}\leq C_F,\,\,\,\,\,\,\|\tilde{F}_0\|_{L^\infty(X_0^\alpha, L_g^2(0,1))}\leq C_F.$$

\item[(d)]$\tilde{F}_\varepsilon$, $\tilde{F}^\eps_0$ and $\tilde{F}_0$ have an uniform bounded support in $\varepsilon\geq 0$, that is:
$$Supp \tilde{F}_\varepsilon\subset \{u\in X_\eps^\alpha :\|u\|_{X_\eps^\alpha}< 2R\},$$
$$Supp \tilde{F}^\eps_0\subset \{u\in X_0^\alpha :\|Eu\|_{X_\eps^\alpha}< 2R\},$$
$$Supp \tilde{F}_0\subset \{u\in X_0^\alpha :\|u\|_{X_0^\alpha}< 2R\},$$
\item[(e)]For all $u\in X_0^\alpha$, 
\begin{equation}\label{rho-beta-*eps}
E\tilde{F}^\eps_0(u)=\tilde{F}_\varepsilon(Eu),\qquad\textrm{and}\quad ED\tilde{F}^\eps_0(u)= D\tilde{F}_\varepsilon(Eu)E.
\end{equation}
and, for any compact set $K\subset X_0^\alpha$, we have, 
\begin{equation}\label{rho-beta-eps0}
\sup_{u_0\in K}\|\tilde{F}_\eps(Eu_0)\mathord-E\tilde{F}_0(u_0)\|_{X_\eps^\alpha}\rightarrow 0,
\end{equation}
\begin{equation}\label{rho-beta-*0}
\sup_{u_0\in K}\|\tilde{F}^\eps_0(u_0)\mathord-\tilde{F}_0(u_0)\|_{X_0^\alpha}\rightarrow 0,
\end{equation}
as $\eps\rightarrow 0$.
\end{itemize}
\end{lem}

\begin{re}
In particular, hypothesis {\bf(H2')} from Section \ref{previous} holds for the three nonlinearities, $\tilde F_\eps$, $\tilde F_0^\eps$ and $\tilde F_0$.  
Moreover, the value of $\rho(\eps)$ and $\beta(\eps)$ from {\bf(H2')}, which depend on the nonlinearities we are considering,  are the following:
$$\rho(\eps), \beta(\eps)=\left\{
\begin{array}{ll}
0& \hbox{ with the nonlinearities } \tilde F_\eps \hbox{ and } \tilde F_0^\eps \\
o(1) & \hbox{ with the nonlinearities } \tilde F_\eps \hbox{ and } \tilde F_0\\
o(1)& \hbox{ with the nonlinearities } \tilde F_0^\eps \hbox{ and } \tilde F_0
\end{array}
\right.
$$
\end{re}

\par\bigskip
\paragraph{\sl Proof. }  (a) This follows directly from definition of $\tilde{F}_\eps$, $\tilde{F}^\eps_0$ and $\tilde{F}_0$, see (\ref{definition-cutoff-eps})-(\ref{definition-cutoff-0}). 

\par\noindent (b) We proceed as follows. Since $F_\eps$ and $\Theta_\eps$ are globally Lipschitz from $X_\eps^\alpha$ to $L^2(Q)$, $\eps>0$, and from $X_0^\alpha$ to $L^2_g(0, 1)$  see Lemma \ref{nonlinearity-C1} and \cite{JamesRobinson}, Lemma 15.7, then $F_\eps$, $\tilde{F}^\eps_0$, $\tilde{F}_0$, are globally Lipschitz from $X_\eps^\alpha$ to $L^2(Q)$ and from $X_0^\alpha$ to $L^2_g(0, 1)$, respectively. So, it remains to prove estimate \ref{Holder}. 

Note that, $D\tilde{F}_\varepsilon(u)=\Theta_\varepsilon(u)DF_\varepsilon(u)+F_\varepsilon(u)D\Theta_\varepsilon(u) $. Then, we can decompose $\|D\tilde{F}_\eps(u)-D\tilde{F}_\eps(v)\|_{\mathcal{L}(X_\eps^\alpha, L^2(Q))}$ as follows,
$$\|D\tilde{F}_\eps(u)-D\tilde{F}_\eps(v)\|_{\mathcal{L}(X_\eps^\alpha, L^2(Q))}\leq$$
$$\|[\Theta_\eps(u)-\Theta_\eps(v)]DF_\eps(u)\|_{\mathcal{L}(X_\eps^\alpha, L^2(Q))}+\|\Theta_\eps(v)[DF_\eps(u)-DF_\eps(v)]\|_{\mathcal{L}(X_\eps^\alpha, L^2(Q))}+$$
$$+\|[F_\eps(u)-F_\eps(v)]D\Theta_\eps(u)\|_{\mathcal{L}(X_\eps^\alpha, L^2(Q))}+\|F_\eps(v)[D\Theta_\eps(u)-D\Theta_\eps(v)]\|_{\mathcal{L}(X_\eps^\alpha, L^2(Q))}=$$
$$=I_1+I_2+I_3+I_4.$$
Since $\Theta_\eps$ is globally Lipschitz with uniform Lipschitz constant, that we denote by  $L_{\hat{\Theta}}$, see \cite[Lemma 15.7]{JamesRobinson} and $\|DF_\eps(u)\|_{\mathcal{L}(X_\eps^\alpha, L^2(Q))}\leq L_F$, see Lemma \ref{nonlinearity-C1}, then
$$I_1\leq L_{\hat{\Theta}} L_F\|u-v\|_{X_\eps^\alpha}. $$
Moreover, by Lemma \ref{nonlinearity-C1} $F_\eps\in C^{1, \theta_F}(X_\eps^\alpha, L^2(Q))$. Hence,
$$I_2\leq L_F\|u-v\|_{X_\eps^\alpha}^{\theta_F}, \qquad\,\,\,\,\textrm{and}\qquad\,\,\,\, I_3\leq L_FL_{\hat{\Theta}}\|u-v\|_{X_\eps^\alpha}.$$
To obtain an estimate for $I_4$, we first calculate the expression for $D\Theta_\eps(u)$. By definition of $\Theta_\eps$, see (\ref{funcionthetaR}), we have for any $u\in X_\eps^\alpha$,
$$D\Theta_\eps(u)=\hat{\Theta}'(\|u\|_{X_\eps^\alpha}^2)2(u, \cdot)_{X_\eps^\alpha},$$
where the function $\hat\Theta$ is defined in \eqref{funcionthetaR},  $'$ is the usual derivative and $(\cdot, \cdot)_{X_\eps^\alpha}$ is the scalar product in the Hilbert space $X^\alpha_\eps$. Hence, 
$$I_4\leq C_F\sup_{\|\varphi\|_{X^\alpha_\eps}=1}\Big\{\Big|\hat{\Theta}'(\|u\|_{X_\eps^\alpha}^2)2(u, \varphi)_{X_\eps^\alpha} - \hat{\Theta}'(\|v\|_{X_\eps^\alpha}^2)2(v, \varphi)_{X_\eps^\alpha}\Big|\Big\}$$
where $C_F$ is the bound from Lemma \ref{nonlinearity-C1} i).  But, 

$$\Big|\hat{\Theta}'(\|u\|_{X_\eps^\alpha}^2)2(u, \varphi)_{X_\eps^\alpha} - \hat{\Theta}'(\|v\|_{X_\eps^\alpha}^2)2(v, \varphi)_{X_\eps^\alpha}\Big|\leq  $$
$$\leq \left|\left(\hat{\Theta}'(\|u\|_{X_\eps^\alpha}^2)-\hat{\Theta}'(\|v\|_{X_\eps^\alpha}^2)\right)2(u, \varphi)_{X_\eps^\alpha}\right|+\left|\Theta'(\|v\|_{X_\eps^\alpha}^2)2(u-v, \varphi)\right|=$$
$$=I_{41}+I_{42}.$$
We first analyze $I_{41}$. Since $\hat{\Theta}$ is a $C^\infty$ function with bounded support in $\mathbb{R}$, then $\hat{\Theta}'$ is globally Lipschitz with Lipschitz constant $L_{\hat{\Theta}}$.
So, 
$$I_{41}\leq 2L_{\hat{\Theta}}\|u\|_{X_\eps^\alpha}\|\varphi\|_{X_\eps^\alpha}\left|\|u\|^2_{X_\eps^\alpha}-\|v\|^2_{X_\eps^\alpha}\right|=$$
$$=2L_{\hat{\Theta}}\|u\|_{X_\eps^\alpha}\left|(\|u\|_{X_\eps^\alpha}+\|v\|_{X_\eps^\alpha})(\|u\|_{X_\eps^\alpha}-\|v\|_{X_\eps^\alpha})\right|\leq$$
$$\leq 2L_{\hat{\Theta}}\|u\|_{X_\eps^\alpha}\left(\|u\|_{X_\eps^\alpha}+\|v\|_{X_\eps^\alpha}\right)\|u-v\|_{X_\eps^\alpha}.$$
We distinguish the following cases:
\begin{itemize}
\item[(1)] If $\|u\|^2_{X_\eps^\alpha}, \|v\|^2_{X_\eps^\alpha}\leq 8R^2$, then
$$I_{41}\leq 32 L_{\hat{\Theta}} R^2\|u-v\|_{X_\eps^\alpha}.$$
\item[(2)] If $\|u\|^2_{X_\eps^\alpha}, \|v\|^2_{X_\eps^\alpha}\geq 8R^2$, then $I_{41}=0$, beacause $\Theta'(\|u\|^2_{X_\eps^\alpha})=\Theta'(\|v\|^2_{X_\eps^\alpha})=0$
\item[(3)] If $\|u\|^2_{X_\eps^\alpha}\leq 8R^2$ and $\|v\|^2_{X_\eps^\alpha}\geq 8R^2$, then we always have $\Theta'(\|v\|^2_{X_\eps^\alpha})=0$. We also distinguish two cases,
\begin{itemize}
\item[(3.1)] If $\|u\|^2_{X_\eps^\alpha}\geq 4R^2$, then again $\Theta'(\|u\|^2_{X_\eps^\alpha})=0$ and therefore $I_{41}=0$.
\item[(3.2)] If $\|u\|^2_{X_\eps^\alpha}\leq 4R^2$, then $\|u-v\|_{X_\eps^\alpha}\geq |\|u\|_{X_\eps^\alpha}-\|v\|_{X_\eps^\alpha}|\geq \frac{1}{2}R$. So, $1\leq \frac{2}{R}\|u-v\|_{X_\eps^\alpha}$, and
$$I_{41}\leq 8R^2|\Theta'(\|u\|_{X_\eps^\alpha}^2)|\leq  16R^2 L_{\hat{\Theta}}\frac{\|u-v\|_{X_\eps^\alpha}}{R}=16RL_{\hat{\Theta}}\|u-v\|_{X_\eps^\alpha}$$

\end{itemize}
\end{itemize}
Therefore, 
$$I_{41}\leq 32L_{\hat{\Theta}} R^2\|u-v\|_{X_\eps^\alpha}.$$
\medskip
Term $I_{42}$ can be directly estimated as follows,
$$I_{42}\leq 2L_{\hat{\Theta}}\|u-v\|_{X_\eps^\alpha}.$$
So
$$I_4\leq (32R^2+2)L_{\hat{\Theta}} \|u-v\|_{X_\eps^\alpha}. $$

Hence, putting all the information together, we get 
$$\|D\tilde{F}_\eps(u)-D\tilde{F}_\eps(v)\|_{\mathcal{L}(X_\eps^\alpha, L^2(Q))}\leq  L_F\|u-v\|_{X_\eps^\alpha}^{\theta_F},$$
with $L_F>0$ independent of $\eps$, as we wanted to prove.

\bigskip

To obtain the same result for $\tilde{F}^\eps_0$ and $\tilde{F}_0$, the proof is exactly the same, step by step.
\par\medskip \noindent (c) This property follows from Lemma \ref{nonlinearity-C1}, item (i).
\par\medskip \noindent (d) It follows directly from the definition of $\Theta_\eps$ and $\Theta_0$.

\par\medskip \noindent (e)
Finally, note that $F_0(u)=f(u(x))=F_\eps(Eu)$. Then, for $u\in X_0^\alpha$,
$$E\tilde{F}^\eps_0(u)=\Theta_\eps^R(Eu)EF_0(u)=\Theta_\eps^R(Eu)f(u(x))=\Theta_\eps^R(Eu)F_\eps(Eu)=\tilde{F}_\eps(Eu),$$
and, since $D\tilde{F}^\eps_0(u)=\Theta_\eps^R(Eu)DF_0(u)+F_0(u)D\Theta^R_\eps(Eu)$,
$$ED\tilde{F}^\eps_0(u)= \Theta_\eps^R(Eu)EDF_0(u)+EF_0(u)D\Theta^R_\eps(Eu)=$$
$$=\Theta_\eps^R(Eu)DF_\eps(Eu)E+F_\eps(Eu)D\Theta^R_\eps(Eu)=D\tilde{F}_\eps(Eu)E.$$
Moreover, for any $u_0\in K\subset X_0^\alpha$ with $K$ compact, we have,


$$\|\tilde{F}_\eps(Eu_0)-E\tilde{F}_0(u_0)\|_{X_\eps}\leq$$
$$\|[\Theta^R_\eps(Eu_0)-\Theta^R_0(u_0)]F_\eps(Eu_0)\|_{X_\eps}+\|\Theta^R_0(u_0)[F_\eps(Eu_0)-EF_0(u_0)]\|_{X_\eps}=$$
$$\|[\Theta^R_\eps(Eu_0)-\Theta^R_0(u_0)]F_\eps(Eu_0)\|_{X_\eps}\leq C_F L_{\hat\Theta}|\|Eu_0\|^2_{X_\eps^\alpha}-\|u_0\|^2_{X_0^\alpha}|=$$
$$C_F L_{\hat\Theta}|( \|Eu_0\|_{X_\eps^\alpha}+\|u_0\|_{X_0^\alpha})(\|Eu_0\|_{X_\eps^\alpha}-\|u_0\|_{X_0^\alpha})|\leq $$
$$C_F L_{\hat\Theta}(2e^2+1)\|u_0\|_{X_0^\alpha}|\|Eu_0\|_{X_\eps^\alpha}-\|u_0\|_{X_0^\alpha}|,$$
in the last inequality we have applied the bound for operator $E$ obtained in \eqref{Ealpha}.

Hence, since $K$ is a compact subset of $X_0^\alpha$, by Lemma \ref{normaproyeccionextension} item (ii),
$$\sup_{u_0\in K}\|\tilde{F}_\eps(Eu_0)-E\tilde{F}_0(u_0)\|_{X_\eps}\rightarrow 0,$$
when $\eps$ tends to zero.

\medskip 
We omit the proof of \eqref{rho-beta-*0} for being equal to the proof of \eqref{rho-beta-eps0}.  \begin{flushright}$\blacksquare$\end{flushright}

\section{Inertial manifolds and reduced systems}\label{InertialManifoldsConstruction}

We present the construction of inertial manifolds for problems (\ref{system-eps}), (\ref{system-*}) and \eqref{system-0}. Remember that, with these manifolds, our problem is reduced to a finite dimensional system.

The existence of these manifolds is guaranteed by the existence of spectral gaps, large enough, in the spectrum of the associated linear elliptic operators, see \cite{Sell&You}. Moreover, these spectral gaps are going to be garanteed by the existence of the spectral gaps for the limiting problem together with the spectral convergence of the linear eliptic operators, which is obtained from {\bf (H1)}, see Section \ref{previous} and \cite{Arrieta-Santamaria-C0,Arrieta-Santamaria-C1}.

With the notations from Section \ref{setting}, by Proposition \ref{resolvente}  for $h_\varepsilon=M f_\varepsilon$ we have, see Remark \ref{version-abstracta}
\begin{equation}\label{resolventefuerte2}
\|A_\varepsilon^{-1}-EA_0^{-1}M\|_{\mathcal{L}(L^2(Q), H^1_{\bm\varepsilon}(Q))}\leq C\varepsilon,
\end{equation}
and for $f_\varepsilon=Eh_\varepsilon$, see Remark \ref{estimaciondebilresolvente},
\begin{equation}\label{resolventedebil}
\|A_\varepsilon^{-1}E-EA_0^{-1}\|_{\mathcal{L}(L^2_g(0, 1), H^1_{\bm\varepsilon}(Q))}\leq C\varepsilon.
\end{equation}

These two estimates imply that hypothesis {\bf(H1)}  holds with $\alpha=1/2$ and therefore it also holds for any $0\leq \alpha\leq 1/2$.  Moreover, the parameter $\tau(\eps)$  is $\tau(\eps)=\eps$.  In the sequel we will use the notation introduce in Section \ref{previous} with respect the eigenvalues, projections, etc..

The limit operator $A_0$ is of Sturm-Liouville type of one dimension. Following \cite{Hale&Raugel3}, Lemma 4.2, we know that there exists $N_0$ such that for all $m\geq N_0$
\begin{equation}\label{eigenvalue-bound}
\pi^2\left(m+\frac{1}{4}\right)^2\leq\lambda^0_m\leq \pi^2\left(m+\frac{3}{4}\right)^2.
\end{equation}
This implies that for $m\geq N_0$,
\begin{equation}\label{eigenvalue-gap}
\pi^2 (m+1)\leq\lambda_{m+1}^0-\lambda_m^0\leq 3\pi^2(m+1).
\end{equation}
%
%

\bigskip
Taking $0<\alpha<1/2$, we get from \eqref{eigenvalue-bound} and \eqref{eigenvalue-gap} that for each $M>0$ large enough, we can choose $m\in \N$ also large enough such that
$$\lambda_{m+1}^0-\lambda_m^0\geq M[(\lambda_{m+1}^0)^\alpha+(\lambda_{m}^0)^\alpha]$$

This means that we are in conditions to apply Proposition \ref{existenciavariedadinercial} 
obtaining that 
%
%
there exist $L<1$ and $0<\varepsilon_1\leq\varepsilon_0$ such that for all $0<\varepsilon\leq\varepsilon_1$ there exist inertial manifolds $\mathcal{M}_\varepsilon$, $\mathcal{M}^\eps_0$ and $\mathcal{M}_0$ for (\ref{system-eps}), (\ref{system-*}) and (\ref{system-0}), given by the ``graph" of functions $\Phi_\varepsilon, \Phi^\eps_0, \Phi_0\in\mathcal{F}_\eps(L, 2R)$, 
\begin{equation}\label{inertialmanifold-eps}
\mathcal{M}_\varepsilon=\{j^{-1}_\varepsilon(z)+\Phi_\varepsilon(z);\,\,z\in\mathbb{R}^m\},
\end{equation}
\begin{equation}\label{inertialmanifold-*}
\mathcal{M}^\eps_0=\{j^{-1}_0(z)+\Phi^\eps_0(z);\,\,z\in\mathbb{R}^m\},
\end{equation}
\begin{equation}\label{inertialmanifold-0}
\mathcal{M}_0=\{j^{-1}_0(z)+\Phi_0(z);\,\,z\in\mathbb{R}^m\},
\end{equation}


\bigskip

If we denote by $T_{\mathcal{M}_\varepsilon}$, $T_{\mathcal{M}^\eps_0}$ and $T_{\mathcal{M}_0}$ the time one maps of the semigroup restricted to the inertial manifolds $\mathcal{M}_\varepsilon$, $\mathcal{M}^\eps_0$ and $\mathcal{M}_0$, respectively, for $u_\varepsilon\in\mathcal{M}_\varepsilon$, $u^\eps_0\in\mathcal{M}^\eps_0$ and $u_0\in\mathcal{M}_0$ and $z\in\mathbb{R}^m$, the time one maps satisfy the following equalities,
$$T_{\mathcal{M}_\varepsilon}(u_\varepsilon)= p_\varepsilon(1)+\Phi_\varepsilon(j_\varepsilon(p_\varepsilon(1))),$$
$$T_{\mathcal{M}^\eps_0}(u^\eps_0)= p^\eps_0(1)+\Phi^\eps_0(j_0(p^\eps_0(1))),$$
$$T_{\mathcal{M}_0}(u_0)= p_0(1)+\Phi_0(j_0(p_0(1))),$$
with $p_\eps(t)$, $p^\eps_0(t)$ and $p_0(t)$ the solutions of
\begin{equation}\label{equationp-eps-modified}
\left\{
\begin{array}{l}
p_t=-A_\varepsilon p+\mathbf{P}_{\mathbf{m}}^{\bm\eps} \tilde{F}_\eps(p+\Phi_\varepsilon( j_\eps (p(t))))\\
p(0)=j_\eps^{-1}(z),
\end{array}
\right.
\end{equation}
\begin{equation}\label{equationp-*-modified}
\left\{
\begin{array}{l}
p_t=-A_0 p+\mathbf{P}_{\mathbf{m}}^{\mathbf{0}} \tilde{F}^\eps_0(p+\Phi^\eps_0( j_0 (p(t))))\\
p(0)=j_0^{-1}(z),
\end{array}
\right.
\end{equation}
\begin{equation}\label{equationp-0-modified}
\left\{
\begin{array}{l}
p_t=-A_0 p+\mathbf{P}_{\mathbf{m}}^{\mathbf{0}} \tilde{F}_0(p+\Phi_0( j_0 (p(t))))\\
p(0)=j_0^{-1}(z).
\end{array}
\right.
\end{equation}

Moreover, $j_\varepsilon(p_\eps(t))$, $j_0(p^\eps_0(t))$ and $j_0(p_0(t))$ satisfy the following systems in $\mathbb{R}^m$,
\begin{equation}\label{equationonRm-eps}
\left\{
\begin{array}{l}
z_t=-j_\varepsilon A_\varepsilon j_\varepsilon^{-1}z+j_\varepsilon\mathbf{P}_{\mathbf{m}}^{\bm\eps} \tilde{F}_\eps(j_\varepsilon^{-1}(z)+\Phi_\varepsilon(z))\\
z(0)=z^0,
\end{array}
\right.
\end{equation}
\begin{equation}\label{equationonRm-*}
\left\{
\begin{array}{l}
z_t=-j_0 A_0 j_0^{-1}z+j_0\mathbf{P}_{\mathbf{m}}^{\mathbf{0}} \tilde{F}^\eps_0(j_0^{-1}(z)+\Phi^\eps_0(z))\\
z(0)=z^0,
\end{array}
\right.
\end{equation}
\begin{equation}\label{equationonRm-0}
\left\{
\begin{array}{l}
z_t=-j_0 A_0 j_0^{-1}z+j_0\mathbf{P}_{\mathbf{m}}^{\mathbf{0}} \tilde{F}_0(j_0^{-1}(z)+\Phi_0(z))\\
z(0)=z^0.
\end{array}
\right.
\end{equation}

We write them in the following way:

\begin{equation}\label{equationonRm-epsH}
\left\{
\begin{array}{l}
z_t=-j_\varepsilon A_\varepsilon j_\varepsilon^{-1}z+H_\eps(z)\\
z(0)=z^0,
\end{array}
\right.
\end{equation}
\begin{equation}\label{equationonRm-*H}
\left\{
\begin{array}{l}
z_t=-j_0 A_0 j_0^{-1}z+H^\eps_0(z)\\
z(0)=z^0,
\end{array}
\right.
\end{equation}
\begin{equation}\label{equationonRm-0H}
\left\{
\begin{array}{l}
z_t=-j_0 A_0 j_0^{-1}z+H_0(z)\\
z(0)=z^0,
\end{array}
\right.
\end{equation}

\bigskip
where
$$H_\eps, H_0^\eps, H_0:\mathbb{R}^m\longrightarrow \mathbb{R}^m,$$
are given by 
$$H_\eps=j_\eps\mathbf{P}_{\mathbf{m}}^{\bm \eps} \tilde{F}_\eps(j_\eps^{-1}(z)+\Phi_\eps(z)), $$
$$H_0^\eps=j_0\mathbf{P}_{\mathbf{m}}^{\mathbf{0}} \tilde{F}^\eps_0(j_0^{-1}(z)+\Phi^\eps_0(z)),$$
$$H_0=j_0\mathbf{P}_{\mathbf{m}}^{\mathbf{0}} \tilde{F}_0(j_0^{-1}(z)+\Phi_0(z)).$$
They are of compact support, 
$$supp(H_\eps), supp(H_0^\eps), supp(H_0)\subset B_{R'},$$ 
and $B_{R'}$ denotes a ball in $\mathbb{R}^m$ of some radius $R'>0$ centered at the origin. 


\bigskip

We have the following result
\begin{prop}\label{MS-Rm}
If all equilibria of \eqref{equationon(01)} are hyperbolic, then the time one map of \eqref{equationonRm-0H} is a Morse-Smale (gradient like) map.
\end{prop}

\paragraph{\sl Proof. }
Since all the equilibrium points of (\ref{equationon(01)}) are hyperbolic, by \cite{Henry2} the stable and unstable manifolds intersect transversally and so, the time one map of the dynamical system generated by (\ref{system-0}) is a Morse-Smale (gradient like) map. In \cite{PilyuginShaDyn}, Section 3.4, S. Y. Pilyugin proves that, then, the time one map $T_{\mathcal{M}_0}$ corresponding to the limit system in the inertial manifold $\mathcal{M}_0$ is a Morse-Smale (gradient like) map in a neighborhood $V$ of the attractor $\mathcal{A}_0$ in this inertial manifold, $V\subset\mathcal{M}_0$. Then, the time one map $\bar{T}_0$ of the limit system in $\mathbb{R}^m$ generated by (\ref{equationonRm-0H}) is Morse-Smale (gradient like).
\begin{flushright}$\blacksquare$\end{flushright}


\section{Rate of the distance of attractors}\label{distanceattractors}
In this section we give an estimate for the distance of attractors related to (\ref{equationonQepsilon}) and (\ref{equationon(01)}), proving our main result, Theorem \ref{maintheorem-thindomain}. To accomplish this, we start showing the following important results about the relation of the time one maps of the dynamical systems related to  (\ref{equationonRm-epsH}), (\ref{equationonRm-*H}) and (\ref{equationonRm-0H}) and the ones corresponding to (\ref{equationon(01)}) and (\ref{equationonQ}).

\bigskip

Let us denote by $\bar{T}_\varepsilon, \bar{T}^\eps_0, \bar{T}_0:\mathbb{R}^m\rightarrow\mathbb{R}^m,$ the time one maps of the dynamical systems generated by (\ref{equationonRm-epsH}), (\ref{equationonRm-*H}) and (\ref{equationonRm-0H}), respectively.

\bigskip

In the following result, we analyze its convergence of these time one maps. 

\begin{lem}\label{distanciasemigruposRm}
We have,
$$\|\bar{T}_\varepsilon-\bar{T}^\eps_0\|_{C^1(\mathbb{R}^m, \mathbb{R}^m)}\rightarrow 0, $$
$$\|\bar{T}^\eps_0-\bar{T}_0\|_{C^1(\mathbb{R}^m, \mathbb{R}^m)}\rightarrow 0, $$
as $\eps\rightarrow 0$.
Moreover, we have,
\begin{equation}\label{estimate-convergence-C0}
\| \bar{T}_\varepsilon-\bar{T}^\eps_0\|_{L^\infty(\mathbb{R}^m, \mathbb{R}^m)}\leq C\varepsilon|\log(\varepsilon)|,
\end{equation}
with $C$ independent of $\varepsilon$.
\end{lem}
\paragraph{\sl Proof. }
Note that $\tilde{F}_\eps\in C^{1, \theta_F}(X_\eps^\alpha, L^2(Q))$, $\tilde{F}_0^\eps, \tilde{F}_0\in C^{1, \theta_F}(X_0^\alpha, L^2_g(0, 1))$, see Lemma \ref{nonlinearity} item (b), and $\Phi_\eps\in C^{1, \theta}(\mathbb{R}^m, X_\eps^\alpha)$, $\Phi_0^\eps, \Phi_0\in C^{1, \theta}(\mathbb{R}^m, X_0^\alpha)$ for certain small $\theta$, see Proposition \ref{FixedPoint-E^1Theta}. Then, it is easy to show that $H_\eps, H_0^\eps, H_0\in C^{1, \theta}(\mathbb{R}^m, \mathbb{R}^m)$ for $\theta>0$ small enough and
\begin{equation}\label{bound-C1theta}
\|H_\eps\|_{C^{1, \theta}(\mathbb{R}^m, \mathbb{R}^m)}, \|H_0^\eps\|_{C^{1, \theta}(\mathbb{R}^m, \mathbb{R}^m)}, \|H_0\|_{C^{1, \theta}(\mathbb{R}^m, \mathbb{R}^m)}\leq \mathbf{M},
\end{equation}
with $\mathbf{M}$ independent of $\eps$.
Moreover, by Lemma \ref{nonlinearity} item (e) we have that, $$\|\tilde{F}_\eps E- E\tilde{F}_0^\eps\|_{C^0(X_0^\alpha, X_\eps)}=0,$$ and for $K=\{u_0=p_0+\Phi_0(p_0)\,\,\, \textrm{with}\,\,\, p_0\in[\varphi_1^0, ..., \varphi_m^0]\,\,\,\,\,\textrm{and}\,\,\,\|p_0\|_{X_0^\alpha}\leq 2R\}\subset X_0^\alpha $
$$\sup_{u_0\in K}\|\tilde{F}_0^\eps(u_0) - \tilde{F}_0(u_0)\|_{X_0}\rightarrow 0,$$
as $\eps\rightarrow 0$. Then, since we have $j_\eps\rightarrow j_0$ and $\mathbf{P}_{\mathbf{m}}^{\bm\eps} \rightarrow \mathbf{P}_{\mathbf{m}}^{\mathbf{0}} $, see 
 Remark 3.3  and Lemma 3.7  and Lemma 5.4 from \cite{Arrieta-Santamaria-C0}, we have that
\begin{equation}\label{convergence-C0}
\|H_\eps-H_0^\eps\|_{C^0(\mathbb{R}^m, \mathbb{R}^m)}\rightarrow 0, \,\,\,\,\,\,\,\,\,\,\,\,\,\,\,\|H_0^\eps-H_0\|_{C^0(\mathbb{R}^m, \mathbb{R}^m)}\rightarrow 0. 
\end{equation}
Hence, \eqref{bound-C1theta}, \eqref{convergence-C0} and the fact that the support is contained in $B_{R'}$ imply 
$$\|H_\eps-H_0^\eps\|_{C^{1, \theta'}(B_{R'}, \mathbb{R}^m)}\rightarrow 0, \,\,\,\,\,\,\,\,\,\,\,\,\,\,\,\|H_0^\eps-H_0\|_{C^{1, \theta'}(B_{R'}, \mathbb{R}^m)}\rightarrow 0, $$
as $\eps\rightarrow 0$, for $\theta'<\theta$. For this, we are using the compact embedding $C^{1, \theta}(B, \mathbb{R}^m)\hookrightarrow C^{1, \theta'}(B, \mathbb{R}^m)$ for all $\theta'<\theta$, the convergence \eqref{convergence-C0} and the boundness of $H_\eps, H_0^\eps, H_0$ in  $C^{1, \theta}(B, \mathbb{R}^m)$. In particular, we have this convergence in the $C^1$-topology.

With this, we obtain the desired convergence,
$$\|\bar{T}_\varepsilon-\bar{T}^\eps_0\|_{C^1(\mathbb{R}^m, \mathbb{R}^m)}\rightarrow 0, $$
$$\|\bar{T}^\eps_0-\bar{T}_0\|_{C^1(\mathbb{R}^m, \mathbb{R}^m)}\rightarrow 0. $$


Now, since systems \eqref{system-eps} and \eqref{system-*} satisfy hypotheses {\bf (H1)} and {\bf (H2')}, then, we can apply the results  from Section \ref{previous}  to obtain estimate \eqref{estimate-convergence-C0}. Hence,

$$\|\bar{T}^\eps_0-\bar{T}_\eps\|_{L^\infty(\mathbb{R}^m, \mathbb{R}^m)}=\sup_{z\in\mathbb{R}^m}|\bar{T}^\eps_0(z)-\bar{T}_\varepsilon(z)|_{0,\alpha}=$$
$$\sup_{z\in\mathbb{R}^m}|z^\eps_0(1)-z_\varepsilon(1)|_{0,\alpha}=\sup_{z\in\mathbb{R}^m}|j_0(p^\eps_0(1))-j_\varepsilon(p_\varepsilon(1))|_{0,\alpha},$$
where $p_\varepsilon(t)$ and $p^\eps_0(t)$ are the solutions of (\ref{equationp-eps-modified}) and (\ref{equationp-*-modified}) with $p_\eps(0)=j_\eps^{-1}(z)$, $p^\eps_0(0)=j_0^{-1}(z)$,  and $z_\eps$, $z^\eps_0$, the solutions of (\ref{equationonRm-eps}) and (\ref{equationonRm-*}) with $z_\eps(0)=z$, $z^\eps_0(0)=z$. 

By Lemma 5.4 from \cite{Arrieta-Santamaria-C0}, and since $\kappa=1$, we obtain,
$$|j_0(p^\eps_0(1))-j_\varepsilon(p_\varepsilon(1))|_{0,\alpha}\leq 2\|Ep^\eps_0(1)-p_\varepsilon(1)\|_{X_\eps^\alpha}+2C_P\varepsilon\|p^\eps_0(1)\|_{L^2_g(0, 1)},$$
with $C_P\sim (\lambda_m^0)^3$ a constant from the estimate of the distance of spectral projections, $\|E\mathbf{P}_{\mathbf{m}}^{\mathbf{0}}-\mathbf{P}_{\mathbf{m}}^{\bm\varepsilon}E\|_{\mathcal{L}(L^2_g(0, 1), X_\eps^\alpha)}$, see Lemma 3.7 from \cite{Arrieta-Santamaria-C0}.

Moreover, since $\|E\tilde{F}^\eps_0- \tilde{F}_\eps E\|_{L^\infty(X_0^\alpha, L^2(Q))}=0$, (see Lemma \ref{nonlinearity}, item (e)) applying Lemma 5.6 from \cite{Arrieta-Santamaria-C0} with $t=1$  and Proposition \ref{resolvente} we have
$$\|Ep^\eps_0(1)-p_\varepsilon(1)\|_{X_\eps^\alpha}\leq C(\|E\Phi^\eps_0-\Phi_\varepsilon\|_{L^\infty(\mathbb{R}^m, X_\eps^\alpha)}+\varepsilon).$$
Then, 
$$\|\bar{T}^\eps_0-\bar{T}_\varepsilon\|_{L^\infty(\mathbb{R}^m, \mathbb{R}^m)}=\sup_{z\in\mathbb{R}^m}|\bar{T}^\eps_0(z)-\bar{T}_\varepsilon(z)|_{0,\alpha}\leq$$
\begin{equation}\label{Tconvergence}
\leq C(\|E\Phi^\eps_0-\Phi_\varepsilon\|_{L^\infty(\mathbb{R}^m, X_\eps^\alpha)}+\varepsilon)\leq C\eps|\log(\eps)|,
\end{equation}
with $C>0$ independent of $\eps$. Last inequality is obtained applying the result on the distance of the inertial manifolds from Section \ref{previous} (see  Theorem \ref{distaciavariedadesinerciales})  

\begin{flushright}$\blacksquare$\end{flushright}

\bigskip

\begin{re}\label{no-rate-convergence}
Note that an estimate for the rate of convergence of $\|\bar{T}_0-\bar{T}_\eps\|_{L^\infty(\mathbb{R}^m, \mathbb{R}^m)}$ and $\|\bar{T}^\eps_0-\bar{T}_0\|_{L^\infty(\mathbb{R}^m, \mathbb{R}^m)}$ is not obtained in a straightforward way. More precisely, the difficulty lies in analyzing the rate of convergence of $\|Eu_0\|_{X_\eps^\alpha}\rightarrow\|u_0\|_{X^\alpha_0}$, see Lemma \ref{normaproyeccionextension} ii).
\end{re}
\bigskip

We now give an estimate for the distance of the time one maps of the dynamical systems generated by (\ref{equationon(01)}) and (\ref{equationonQ})
\begin{lem}\label{NonLinearSemigroup}
Let $T_0$ and $T_\eps$, $0<\eps\leq\eps_0$, the time one maps corresponding to (\ref{equationon(01)}) and (\ref{equationonQ}), respectively. Then, for $R>0$ large enough, there exists a constant $C=C(R)$ such that for any $w_0\in L^2_g(0,1)$, with $\|w_0\|_{L^2_g(0,1)}\leq R$, we have,
$$\|T_\eps(Ew_0)-ET_0(w_0)\|_{H^1_{\bm\eps}(Q)}\leq C\eps|\log(\eps)|.$$

\end{lem}

\paragraph{\sl Proof. } We have denoted previously by $S_\eps(t)$ and $S_0(t)$ the nonlinear semigroups generated by (\ref{equationonQ}) and (\ref{equationon(01)}) respectively, so that $T_\eps=S_\eps(1)$ and $T_0=S_0(1)$. Hence, with the variation of constants formula,  for $0<t\leq 1$,

$$\| S_\eps(t)(Ew_0)-E S_0(t)(w_0)\|_{H^1_{\bm\eps}(Q)}\leq\|(e^{-A_\eps t}E-Ee^{-A_0 t})w_0\|_{H^1_{\bm\eps}(Q)}+$$
$$+\int_0^t \left\|e^{-A_\eps(t-s)} {F}_\eps(S_\eps(s)Ew_0)-Ee^{-A_0(t-s)} {F}_0(S_0(s)w_0)\right\|_{H^1_{\bm\eps}(Q)}ds\leq$$
$$\leq \|(e^{-A_\eps t}E-Ee^{-A_0 t})w_0\|_{H^1_{\bm\eps}(Q)}+$$
$$+\int_0^t \left\|\left(e^{-A_\eps(t-s)}E- Ee^{-A_0(t-s)}\right) F_0(S_0(s)w_0)\right\|_{H^1_{\bm\eps}(Q)}ds+$$
$$+\int_0^t \left\|e^{-A_\eps(t-s)}\left( F_\eps(ES_0(s)w_0)-F_0(S_0(s)w_0)\right)\right\|_{H^1_{\bm\eps}(Q)}ds+$$
$$+ \int_0^t \left\|e^{-A_\eps(t-s)}\left(F_\eps(S_\eps(s)Ew_0)-F_\eps(ES_0(s)w_0)\right)\right\|_{H^1_{\bm\eps}(Q)}ds.$$
 
But notice that since both $F_\eps$ and $F_0$ are Nemitskii operators of the same function $f:\mathbb{R}\to \mathbb{R}$ then $ {F}_\eps(ES_0(s)w_0)={F}_0(S_0(s)w_0)$, and the third term is identically 0. 

Now, since hypothesis {\bf (H1)} is satisfied, applying Lemma 3.9, Lemma 3.10  from \cite{Arrieta-Santamaria-C0}, 
Lemma \ref{normaproyeccionextension}, Proposition \ref{resolvente} and with Gronwall-Henry inequality, see \cite{Henry1} Section 7, for $t =1$, we obtain,
$$\|T_\eps(Ew_0)-ET_0(w_0)\|_{H^1_{\bm\eps}(Q)}=\|S_\eps(1)(Ew_0)-E S_0(1)(w_0)\|_{H^1_{\bm\eps}(Q)}\leq C \eps|\log(\eps)|, $$
with $C>0$ independent of $\eps$.
\begin{flushright}$\blacksquare$\end{flushright}
\par\bigskip

We show the time one maps are Lipschitz from $L^2(Q)$ to $H^1_{\bm\eps}(Q)$ uniformly in $\eps$.
\begin{lem}\label{timeonemap-lipschitz}
There exists a constant $C>0$ independent of $\eps$ so that, for $0\leq\eps\leq \eps_0$,
$$\|T_\eps(u_\eps)-T_\eps(w_\eps)\|_{H^1_{\bm\eps}(Q)}\leq C\|u_\eps-w_\eps\|_{L^2(Q)}.$$
\end{lem}
\paragraph{\sl Proof. }
By the variation of constants formula, for $0<t\leq 1$, we have 
$$\|S_\eps(t)u_\eps-S_\eps(t)w_\eps\|_{H^1_{\bm\eps}(Q)}\leq \| e^{-A_\eps t} (u_\eps-w_\eps)\|_{H^1_{\bm\eps}(Q)}+$$
$$+\int_0^t\left\|e^{-A_\eps(t-s)}(F_\eps(S_\eps(s)u_\eps)- F_\eps(S_\eps(s)w_\eps))\right\|_{H^1_{\bm\eps}(Q)}ds.$$
Applying Lemma 3.1 from \cite{Arrieta-Santamaria-C0} and Lemma \ref{nonlinearity-C1}, item (ii),
$$\|S_\eps(t)u_\eps-S_\eps(t)w_\eps\|_{H^1_{\bm\eps}(Q)}\leq Ce^{-\lambda_1^\eps t}t^{-\frac{1}{2}}\|u_\eps-w_\eps\|_{L^2(Q)}+ $$
$$+CL_F e^{-\lambda_1^\eps t}\int_0^t e^{\lambda_1^\eps s}(t-s)^{-\frac{1}{2}}\|S_\eps(s)u_\eps-S_\eps(s)w_\eps\|_{H^1_{\bm\eps}(Q)} ds.$$

Applying Gronwall inequality,  for $0<t\leq 1$, we have
$$\|S_\eps(t)u_\eps-S_\eps(t)w_\eps\|_{H^1_{\bm\eps}(Q)}\leq  C t^{-\frac{1}{2}}\|u_\eps-w_\eps\|_{L^2(Q)} e^{-\lambda_1^\eps t},$$
with $C>0$ independent of $\eps$.

Then, for the time one map $T_\eps=S_\eps(1)$ we obtain
$$\|T_\eps(u_\eps)-T_\eps(w_\eps)\|_{H^1_{\bm\eps}(Q)}\leq  C \|u_\eps-w_\eps\|_{L^2(Q)},$$
with $C>0$ independent of $\eps$, which shows the result.
\begin{flushright}$\blacksquare$\end{flushright}

\bigskip

We proceed to prove the main result of this work.

\paragraph{Proof of Theorem \ref{maintheorem-thindomain}}
We obtain now a rate of convergence of attractors $\mathcal{A}_0$ and $\mathcal{A}_\varepsilon$ of the dynamical systems generated by (\ref{equationon(01)}) and (\ref{equationonQ}), respectively. We know that for any $u_0\in\mathcal{A}_0$ and any $u_\eps\in\mathcal{A}_\eps$ there exist a $w_0\in\mathcal{A}_0$ and $w_\eps\in\mathcal{A}_\eps$ such that,
$$u_0=T_0(w_0),\qquad\textrm{and}\quad u_\eps=T_\eps(w_\eps),$$
with $T_0$ and $T_\eps$ the time one maps corresponding to (\ref{equationon(01)}) and (\ref{equationonQ}).

Moreover, as we have said before, for each $\eps>0$ the attractor $\mathcal{A}_\eps$ is contained in the inertial manifold $\mathcal{M}_\eps$ and $\mathcal{A}_0$ is contained in the inertial manifolds $\mathcal{M}^\eps_0$ and $\mathcal{M}_0$. We also have that although $\mathcal{M}_\eps$, $\mathcal{M}^\eps_0$ and $\mathcal{M}_0$ are manifolds close enough, we only can provide explicit rates of the distance between $\mathcal{M}_\eps$ and $\mathcal{M}^\eps_0$ as $\eps$ goes to zero. 

The Hausdorff distance of attractors $\mathcal{A}_0$ and $\mathcal{A}_\eps$ in $H^1_{\bm\eps}(Q)$, is given by
$$dist_{H^1_{\bm\eps}(Q)}(\mathcal{A}_0, \mathcal{A}_\varepsilon)=max\{\sup_{u_0\in\mathcal{A}_0}\inf_{u_\varepsilon\in\mathcal{A}_\varepsilon}\|Eu_0-u_\varepsilon\|_{H^1_{\bm\varepsilon}(Q)}, \sup_{u_\varepsilon\in\mathcal{A}_\varepsilon}\inf_{u_0\in\mathcal{A}_0}\|u_\varepsilon-E u_0\|_{H^1_{\bm\varepsilon}(Q)}\}.  $$

Then, we consider $w_\eps\in\mathcal{A}_\eps$, $0<\eps\leq\eps_0$, given by $w_\varepsilon=j_\varepsilon^{-1}(z_\eps) +\Phi_\varepsilon(z_\eps)$ and $w_0\in\mathcal{A}_0$, given by $w_0=j_0^{-1}(z_0) +\Phi^\eps_0(z_0)$   with $z_\eps\in\bar{\mathcal{A}}_\eps$ and $z_0\in\bar{\mathcal{A}}_0,$ the ``projected'' attractors in $\mathbb{R}^m$ corresponding to (\ref{equationonRm-eps}) and (\ref{equationonRm-*}), respectively. 

We know,
$$\|Eu_0-u_\eps\|_{H^1_{\bm\eps}(Q)}=\|ET_0(w_0)-T_\eps(w_\eps)\|_{H^1_{\bm\eps}(Q)}\leq$$
$$\leq\|ET_0(w_0)-T_\eps(Ew_0)\|_{H^1_{\bm\eps}(Q)}+\|T_\eps(Ew_0)-T_\eps(w_\eps)\|_{H^1_{\bm\eps}(Q)}.$$
Applying Lemma \ref{NonLinearSemigroup} and Lemma \ref{timeonemap-lipschitz}, we have
$$\|Eu_0-u_\eps\|_{H^1_{\bm\eps}(Q)}\leq C \eps|\log(\eps)|+ C\|Ew_0-w_\eps\|_{X_\eps^\alpha}.$$

So, we need to estimate the norm $\|Ew_0-w_\varepsilon\|_{X_\eps^\alpha}$,
where,
$$w_\varepsilon=j_\varepsilon^{-1}(z_\varepsilon)+\Phi_\varepsilon(z_\varepsilon),\qquad z_\varepsilon\in\bar{\mathcal{A}}_\varepsilon,$$
and
$$w_0=j_0^{-1}(z_0)+\Phi^\eps_0(z_0),\qquad z_0\in\bar{\mathcal{A}}_0$$
with $\bar{\mathcal{A}}_\varepsilon$ and $\bar{\mathcal{A}}_0$ the attractors corresponding to (\ref{equationonRm-eps}) and (\ref{equationonRm-*}).

\bigskip

Hence, since $j_0^{-1}(z_0)=\sum_{i=1}^mz_i^0\psi_i^0$ and $j_\eps^{-1}(z_\eps)=\sum_{i=1}^mz_i^\eps\psi_i^\eps$,
$$\|E w_0-w_\varepsilon\|_{X_\eps^\alpha}\leq\|Ej_0^{-1}(z_0)-j_\varepsilon^{-1}(z_\varepsilon)\|_{X_\eps^\alpha} + \|E\Phi^\eps_0(z_0)-\Phi_\varepsilon(z_\varepsilon)\|_{X_\eps^\alpha}\leq$$
$$\leq \|\sum_{i=1}^m(z_i^0-z_i^\varepsilon)E\psi_i^0\|_{X_\eps^\alpha}+\|\sum_{i=1}^mz_i^\varepsilon(E\psi_i^0-\psi_i^\varepsilon)\|_{X_\eps^\alpha}+$$
$$+\|E\Phi^\eps_0(z_0)-E\Phi^\eps_0(z_\varepsilon)\|_{X_\eps^\alpha}+ \|E\Phi^\eps_0(z_\varepsilon)-\Phi_\varepsilon(z_\varepsilon)\|_{X_\eps^\alpha}\leq$$
$$\leq 2|z_0-z_\varepsilon|_{0,\alpha} + \sup_{z_\varepsilon\in\bar{\mathcal{A}}_\varepsilon} |z_\varepsilon|\|E\mathbf{P}_{\mathbf{m}}^{\mathbf{0}}-\mathbf{P}_{\mathbf{m}}^{\bm\varepsilon}E\|_{\mathcal{L}(L^2_g(0, 1), X_\eps^\alpha)}+\|E\Phi^\eps_0-\Phi_\varepsilon\|_{L^\infty(\mathbb{R}^m, X_\eps^\alpha)}.$$
In the last inequality we have applied the estimate of the norm of operator $E$, see \eqref{Ealpha}.

Since $z_0\in\bar{\mathcal{A}}_0$ and $z_\varepsilon\in\bar{\mathcal{A}}_\varepsilon$, then
$$\|E w_0-w_\varepsilon\|_{X_\eps^\alpha}\leq 2|z_0-z_\varepsilon|_{0,\alpha}+|z_\varepsilon|\|E\mathbf{P}_{\mathbf{m}}^{\mathbf{0}}-\mathbf{P}_{\mathbf{m}}^{\bm\varepsilon}E\|_{\mathcal{L}(L^2_g(0, 1), X_\eps^\alpha)}+$$
$$+\|E\Phi^\eps_0-\Phi_\varepsilon\|_{L^\infty(\mathbb{R}^m, X_\eps^\alpha)}=I_1+I_2+I_3.$$
To estimate $I_2$, note that we have studied the convergence of $\|E\mathbf{P}_{\mathbf{m}}^{\mathbf{0}}-\mathbf{P}_{\mathbf{m}}^{\bm\varepsilon}E\|_{\mathcal{L}(L^2_g(0, 1), X_\eps^\alpha)}$ in terms of the distance of the resolvent operators, see \eqref{convergence-of-projection} or \cite[Lemma 3.7]{Arrieta-Santamaria-C0}. Then, in our case, we have,
$$I_2\leq C\varepsilon.$$
By Theorem \ref{distaciavariedadesinerciales}, 
$$I_3\leq C\eps|\log(\eps)|.$$

Hence, putting everything together,
$$\|E w_0-w_\varepsilon\|_{X_\eps^\alpha}\leq 4e^2|z_0-z_\varepsilon|_{0,\alpha}+ C \varepsilon|\log(\varepsilon)|,$$
with $C$ independent of $\varepsilon$.
Then,
$$\sup_{w_0\in\mathcal{A}_0}\inf_{ w_\eps\in\mathcal{A}_\eps}\|E w_0-w_\varepsilon\|_{X_\eps^\alpha}\leq 4e^2\sup_{z_0\in\bar{\mathcal{A}}_0}\inf_{ z_\eps\in\bar{\mathcal{A}}_\eps}|z_0-z_\varepsilon|_{0,\alpha}+ C \varepsilon|\log(\varepsilon)|.$$
Hence,
$$dist_{H^1_{\bm\eps}(Q)}(\mathcal{A}_0, \mathcal{A}_\eps)\leq  4e^2 dist_{\mathbb{R}^m}(\bar{\mathcal{A}}_0, \bar{\mathcal{A}}_\eps)+C \varepsilon|\log(\varepsilon)|.$$

To estimate $dist_H(\bar{\mathcal{A}}_0, \bar{\mathcal{A}}_\eps)$, we need to apply techniques of Shadowing Theory described in Appendix \ref{shadowing}. First, we have by Proposition \ref{MS-Rm}, that the time one map of the system given by the ordinary differential equation \eqref{equationonRm-0H} is a Morse-Smale map. Moreover, by Lemma \ref{distanciasemigruposRm}, we can take $\eps$ small enough so that the time one maps corresponding to \eqref{equationonRm-epsH} and \eqref{equationonRm-*H}, $\bar{T}_\eps$ and $\bar{T}_0^\eps$, respectivelly belong to a $C^1$ neighborhood of $\bar{T}_0$. Then, by Corollary \ref{MS-dist} 
$$dist_{\mathbb{R}^m}(\bar{\mathcal{A}}_0, \bar{\mathcal{A}}_\varepsilon)\leq L\|\bar{T}^\eps_0-\bar{T}_\eps\|_{L^\infty(\mathbb{R}^m, \mathbb{R}^m)},$$
with $L>0$ independent of $\eps$. Hence, using the estimate obtained in Lemma \ref{distanciasemigruposRm},
$$dist_{\mathbb{R}^m}(\bar{\mathcal{A}}_0, \bar{\mathcal{A}}_\eps)\leq C\eps|\log(\eps)|,$$

Putting all together, we get 
$$dist_{H^1_{\bm\eps}(Q)}(\mathcal{A}_0, \mathcal{A}_\varepsilon)\leq C \varepsilon|\log(\varepsilon)|,$$

Finally, applying identity \eqref{norm-relations2}, we have,
$$dist_{H^1(Q_\eps)}(\mathcal{A}_0, \mathcal{A}_\varepsilon)=\eps^{\frac{d-1}{2}}dist_H(\mathcal{A}_0, \mathcal{A}_\varepsilon)\leq C \varepsilon^{\frac{d+1}{2}}|\log(\varepsilon)|,$$
with $C$ independent of $\varepsilon$.
This shows Theorem  \ref{maintheorem-thindomain}. 
  \begin{flushright}$\blacksquare$\end{flushright}

\appendix

\section{Appendix: Proof of Proposition \ref{resolvente}}
\label{proof-resolvente}

We provide in this appendix the proof of the estimates of the resolvent operators contained in Proposition \ref{resolvente}. 

\par\bigskip
\paragraph{\sl Proof. }
The proof of this result follows similar ideas as the proof of Proposition A.8 from \cite{dumbel1}.

Remember that
$$Q_\varepsilon =\{(x, \varepsilon\mathbf{y})\in\mathbb{R}^d: (x, \mathbf{y})\in Q\},$$
where
$$Q=\{(x, \mathbf{y})\in\mathbb{R}^d: 0\leq x\leq1,\; \; \mathbf{y}\in\Gamma^1_x\},$$
and
$$H^1_{\bm\varepsilon}(Q):=(H^1(Q), \|\cdot\|_{H^1_{\bm\varepsilon}(Q)}),$$
with the norm
$$\|u\|_{H^1_{\bm\varepsilon}(Q)}:=\left(\int_Q(|\nabla_xu|^2+\frac{1}{\varepsilon^2}|\nabla_{\mathbf{y}}u|^2+|u|^2)dxd\mathbf{y}\right)^{1/2}.$$
So, by the change of variable theorem,
$$\|u\|_{L^2(Q_\varepsilon)}=\varepsilon^{\frac{d-1}{2}}\|\mathbf{i}_{\bm\varepsilon} u\|_{L^2(Q)},$$
and
$$\|u\|_{H^1(Q_\varepsilon)}=\varepsilon^{\frac{d-1}{2}}\|\mathbf{i}_{\bm\varepsilon} u\|_{H^1_{\bm\varepsilon}(Q)}.$$
Hence, proving this Proposition is equivalent to prove the estimate
$$\|w_\varepsilon-E_\varepsilon v_\varepsilon\|_{H^1(Q_\varepsilon)}\leq C\varepsilon\|f_\varepsilon\|_{L^2(Q_\varepsilon)},$$
where $w_\varepsilon$ and $v_\varepsilon$ are the solutions of the following linear problems, respectivelly,
\begin{equation}
\left\{
\begin{array}{c l r}
-\Delta w_\varepsilon+\mu w_\varepsilon\;&\; = f_\varepsilon, \;&\;\textrm{in}\quad Q_\varepsilon\\
\frac{\partial w_\varepsilon}{\partial\nu_\varepsilon}\;&\;=0\quad&\textrm{on}\quad \partial Q_\varepsilon,
\end{array}
\right.
\end{equation} 
and 
\begin{equation}
\left\{
\begin{array}{r l r}
-\frac{1}{g}(g {v_\varepsilon}_x)_x + \mu v_\varepsilon\;&\; = M_\varepsilon f_\varepsilon, \;&\;\textrm{in}\quad (0, 1)\\
{v_\varepsilon}_x(0)\;&\;=0,\;&\;{v_\varepsilon}_x(1)=0,
\end{array}
\right.
\end{equation} 
with $f_\varepsilon\in L^2(Q_\varepsilon)$. Observe that $u_\varepsilon(x, \mathbf{y})=w_\varepsilon(x, \varepsilon \mathbf{y})$.
 
It is known that the minima
\begin{equation}\label{minimoQepsilon}
\lambda_\varepsilon:=\displaystyle\min_{\varphi\in H^1(Q_\varepsilon)}\left\{\frac{1}{2}\int_{Q_\varepsilon} (|\nabla\varphi|^2 + \mu|\varphi|^2)ds-\int_{Q_\varepsilon}f_\varepsilon\varphi ds\right\},
\end{equation}
\begin{equation}\label{minimo01}
\tau_\varepsilon:=\displaystyle\min_{\varphi\in H_g^1(0,1)}\left\{\frac{1}{2}\int_0^1(g|\varphi'|^2+ g\mu|\varphi|^2)dx-\int_0^1gM_\varepsilon f_\varepsilon\varphi dx\right\},
\end{equation}
with $s=(x, \mathbf{y})\in Q_\varepsilon$, are unique and they are attained at the solutions $w_\varepsilon$ and $v_\varepsilon$. We want to compare both solutions $w_\varepsilon$ and $v_\varepsilon$.  We start by taking the function $v_\eps$ as a test function in (\ref{minimoQepsilon}). We have,
$$\lambda_\varepsilon\leq\frac{1}{2}\int_{Q_\varepsilon}(|\nabla v_\varepsilon|^2+\mu|v_\varepsilon|^2)ds - \int_{Q_\varepsilon}f_\varepsilon v_\varepsilon ds=$$
$$=\frac{1}{2}\int_0^1\int_{\Gamma_x^\varepsilon}(|{v_\varepsilon}_x|^2+\mu|v_\varepsilon|^2)d\mathbf{y}dx-\int_0^1\int_{\Gamma_x^\varepsilon}f_\varepsilon d\mathbf{y}v_\varepsilon dx=$$
$$=\frac{1}{2}\int_0^1|\Gamma_x^\varepsilon|(|{v_\varepsilon}_x|^2+\mu|v_\varepsilon|^2)dx-\int_0^1|\Gamma_x^\varepsilon|M_\varepsilon f_\varepsilon(x, \mathbf{y})v_\varepsilon dx=$$
$$=\varepsilon^{d-1}\left(\frac{1}{2}\int_0^1g(x)(|{v_\varepsilon}_x|^2+\mu|v_\varepsilon|^2)dx-\int_0^1g(x)M_\varepsilon f_\varepsilon(x, \mathbf{y})v_\varepsilon dx\right)=\varepsilon^{d-1}\tau_\varepsilon.$$
That is, we have obtained the estimate,
$$\lambda_\varepsilon\leq\varepsilon^{d-1}\tau_\varepsilon.$$
To look for a lower bound we proceed as follows,
$$\lambda_\varepsilon=\frac{1}{2}\int_{Q_\varepsilon}(|\nabla w_\varepsilon|^2+\mu|w_\varepsilon|^2)ds - \int_{Q_\varepsilon}f_\varepsilon w_\varepsilon ds=$$
$$=\frac{1}{2}\int_{Q_\varepsilon}(|\nabla w_\varepsilon-\nabla v_\varepsilon + \nabla v_\varepsilon|^2+\mu |w_\varepsilon-v_\varepsilon + v_\varepsilon|^2)ds - \int_{Q_\varepsilon }f_\varepsilon(w_\varepsilon-v_\varepsilon + v_\varepsilon) ds=$$
$$=\frac{1}{2}\int_{Q_\varepsilon}(|\nabla w_\varepsilon-\nabla v_\varepsilon|^2+|\nabla v_\varepsilon|^2+2(\nabla w_\varepsilon- \nabla v_\varepsilon)\nabla v_\varepsilon)ds +$$
$$+ \frac{1}{2}\int_{Q_\varepsilon}\mu(|w_\varepsilon- v_\varepsilon|^2 + |v_\varepsilon|^2 + 2(w_\varepsilon-v_\varepsilon)v_\varepsilon)ds-\int_{Q_\varepsilon}f_\varepsilon(w_\varepsilon-v_\varepsilon)ds-\int_{Q_\varepsilon}f_\varepsilon v_\varepsilon ds.$$
From above, we know that $\frac{1}{2}\int_{Q_\varepsilon}(|\nabla v_\varepsilon|^2+\mu|v_\varepsilon|^2)ds-\int_{Q_\varepsilon}f_\varepsilon v_\varepsilon ds=\varepsilon^{d-1}\tau_\varepsilon$, then
$$\lambda_\varepsilon=\frac{1}{2}\int_{Q_\varepsilon}(|\nabla w_\varepsilon-\nabla v_\varepsilon|^2+2(\nabla w_\varepsilon- \nabla v_\varepsilon)\nabla v_\varepsilon)ds + \frac{1}{2}\int_{Q_\varepsilon}\mu(|w_\varepsilon- v_\varepsilon|^2 + 2(w_\varepsilon-v_\varepsilon)v_\varepsilon)ds$$
$$-\int_{Q_\varepsilon}f_\varepsilon(w_\varepsilon-v_\varepsilon)ds + \varepsilon^{d-1}\tau_\varepsilon.$$
To analyze this, we write the last equality like this,
$$\lambda_\varepsilon= \frac{1}{2}\int_{Q_\varepsilon}(|\nabla w_\varepsilon-\nabla v_\varepsilon|^2 + \mu|w_\varepsilon -v_\varepsilon|^2)ds+I_1+I_2-I_3+\varepsilon^{d-1}\tau_\varepsilon,$$
with,
$$I_1:=\int_{Q_\varepsilon}(\nabla w_\varepsilon - \nabla v_\varepsilon)\nabla v_\varepsilon ds,\qquad I_2:=\int_{Q_\varepsilon}\mu(w_\varepsilon-v_\varepsilon)v_\varepsilon ds,$$
and
$$I_3:= \int_{Q_\varepsilon}f_\varepsilon(w_\varepsilon-v_\varepsilon)ds.$$
If we analyze each term with detail, we observe the following,
$$I_1=\int_{Q_\varepsilon}(\nabla w_\varepsilon - \nabla v_\varepsilon)\nabla v_\varepsilon ds=\int_{Q_\varepsilon}\left(\frac{\partial w_\varepsilon}{\partial x}- v'_\varepsilon\right)v'_\varepsilon ds =$$
$$\int_{Q_\varepsilon}\left(M_\varepsilon\frac{\partial w_\varepsilon}{\partial x}-v_\varepsilon'\right)v'_\varepsilon ds=\int_{Q_\varepsilon}\left(M_\varepsilon\frac{\partial w_\varepsilon}{\partial x}-(M_\varepsilon w_\varepsilon)'\right)v'_\varepsilon dx+\int_{Q_\varepsilon}\left((M_\varepsilon w_\varepsilon)'-v'_\varepsilon\right)v'_\varepsilon ds=$$
$$=\int_{Q_\varepsilon}\left(M_\varepsilon\frac{\partial w_\varepsilon}{\partial x}-(M_\varepsilon w_\varepsilon)'\right)v'_\varepsilon ds+\varepsilon^{d-1}\int_0^1g(x)(M_\varepsilon w_\varepsilon-v_\varepsilon)'v'_\varepsilon dx.$$
Since $v_\varepsilon=v_\varepsilon(x)$, we have,
$$I_2=\int_{Q_\varepsilon}\mu(w_\varepsilon-v_\varepsilon)v_\varepsilon ds=\varepsilon^{d-1}\int_0^1\mu g(x)(M_\varepsilon w_\varepsilon-v_\varepsilon)v_\varepsilon dx,$$
and
$$I_3=\int_{Q_\varepsilon}(f_\varepsilon-M_\varepsilon f_\varepsilon)(w_\varepsilon-v_\varepsilon)ds+\int_{Q_\varepsilon}M_\varepsilon f_\varepsilon(w_\varepsilon-v_\varepsilon)=$$
$$=\int_{Q_\varepsilon}(f_\varepsilon-M_\varepsilon f_\varepsilon)(w_\varepsilon-v_\varepsilon)ds+\int_{Q_\varepsilon}M_\varepsilon f_\varepsilon(M_\varepsilon w_\varepsilon-v_\varepsilon)ds=$$
$$=\int_{Q_\varepsilon}(f_\varepsilon-M_\varepsilon f_\varepsilon)(w_\varepsilon-v_\varepsilon)ds+ \varepsilon^{d-1}\int_0^1g(x) M_\varepsilon (f_\varepsilon)(M_\varepsilon w_\varepsilon-v_\varepsilon)dx.$$
That is,
$$I_1=\tilde{I}_1+\varepsilon^{d-1}\int_0^1g(x)\left(M_\varepsilon w_\varepsilon-v_\varepsilon\right)'v'_\varepsilon dx,\qquad I_2=\varepsilon^{d-1}\int_0^1 \mu g(x)(M_\varepsilon w_\varepsilon-v_\varepsilon)v_\varepsilon dx,$$
and
$$I_3=  \tilde{I}_3+\varepsilon^{d-1}\int_0^1g(x)(M_\varepsilon w_\varepsilon-v_\varepsilon)M_\varepsilon f_\varepsilon dx,$$
where $$\tilde{I}_1=\int_{Q_\varepsilon}\left(M_\varepsilon\frac{\partial w_\varepsilon}{\partial x}-\left(M_\varepsilon w_\varepsilon\right)'\right)v'_\varepsilon ds,\qquad\textrm{and}\qquad\tilde{I}_3=\int_{Q_\varepsilon}\left(f_\varepsilon-M_\varepsilon f_\varepsilon\right)(w_\varepsilon-v_\varepsilon)ds.$$
We know that,
$$\int_0^1\left[g(x)(M_\varepsilon w_\varepsilon-v_\varepsilon)'v'_\varepsilon+ \mu g(x)(M_\varepsilon w_\varepsilon-v_\varepsilon)v_\varepsilon\right]dx=\int_0^1g(x)(M_\varepsilon w_\varepsilon-v_\varepsilon)M_\varepsilon f_\varepsilon dx,$$
then,
$$I_1 + I_2 - I_3=\tilde{I}_1-\tilde{I}_3.$$
So, we only need to estimate $\tilde{I}_1$ and $\tilde{I}_3$.

We start with $\tilde{I}_1$.
$$\tilde{I}_1=\int_{Q_\varepsilon}\left(M_\varepsilon\frac{\partial w_\varepsilon}{\partial x}-\left(M_\varepsilon w_\varepsilon\right)'\right)v'_\varepsilon ds.$$
Then, we first estimate $M_\varepsilon\frac{\partial w_\varepsilon}{\partial x}-\left(M_\varepsilon w_\varepsilon\right)'$. For that, we study $(M_\varepsilon w_\varepsilon)'$.
$$(M_\varepsilon w_\varepsilon)'=\frac{d}{dx}\left(\frac{1}{|\Gamma_x^\varepsilon|}\int_{\Gamma_x^\varepsilon}w_\varepsilon(x, \mathbf{y})d\mathbf{y}\right),$$
and by the Change of Variable Theorem with $\mathbf{y}=\varepsilon \mathbf{L}_x(z)$, see (\ref{difeomorfismo}), and $z\in B(0,1)$ the unit ball in $\mathbb{R}^{d-1}$, we have
$$\frac{1}{|\Gamma_x^\varepsilon|}\int_{\Gamma_x^\varepsilon}w_\varepsilon(x, \mathbf{y})d\mathbf{y}=\int_{B(0,1)}w_\varepsilon(x, \varepsilon \mathbf{L}_x(z))\frac{J_{\mathbf{L}_x(z)}}{|\Gamma_x^1|} dz,$$
where $J_{\mathbf{L}_x(z)}$ is the Jacobian of $\mathbf{L}_x$. So,
$$ (M_\varepsilon w_\varepsilon)'=\frac{d}{dx}\left(\int_{B(0, 1)}w_\varepsilon(x, \varepsilon \mathbf{L}_x(z))\frac{J_{\mathbf{L}_x(z)}}{|\Gamma_x^1|} dz\right)=$$
$$=\int_{B(0,1)}\frac{\partial w_\varepsilon}{\partial x} (x, \varepsilon \mathbf{L}_x(z))\frac{J_{\mathbf{L}_x(z)}}{|\Gamma_x^1|}dz + \int_{B(0,1)}\nabla_{\mathbf{y}}w_\varepsilon(x, \varepsilon \mathbf{L}_x(z))\varepsilon\frac{\partial}{\partial x}(\mathbf{L}_x(z))\frac{J_{\mathbf{L}_x(z)}}{|\Gamma_x^1|}dz+$$
$$+\int_{B(0,1)}w_\varepsilon(x, \varepsilon \mathbf{L}_x(z))\frac{\partial}{\partial x}\left(\frac{J_{\mathbf{L}_x(z)}}{|\Gamma_x^1|}\right)dz.$$
To estimate the right side of the above equality, we study each integral separately. We begin with the first one. Undoing the change of variable $\mathbf{y}=\varepsilon \mathbf{L}_x(z)$, 
$$\int_{B(0,1)}\frac{\partial w_\varepsilon}{\partial x} (x, \varepsilon \mathbf{L}_x(z))\frac{J_{L_x(z)}}{|\Gamma_x^1|}dz=\frac{1}{|\Gamma^\varepsilon_x|}\int_{\Gamma_x^\varepsilon}\frac{\partial w_\varepsilon}{\partial x}(x, \mathbf{y})d\mathbf{y}= M_\varepsilon\frac{\partial w_\varepsilon}{\partial x}.$$
For the second integral we use $\left|\frac{\partial \mathbf{L}_x(z)}{\partial x}\right|\leq C$, 
$$\left|\int_{B(0,1)}\nabla_{\mathbf{y}}w_\varepsilon(x, \varepsilon \mathbf{L}_x(z))\varepsilon\frac{\partial}{\partial x}(\mathbf{L}_x(z))\frac{J_{\mathbf{L}_x(z)}}{|\Gamma_x^1|}dz\right|\leq C\varepsilon\int_{B(0,1)}|\nabla_{\mathbf{y}}w_\varepsilon(x, \varepsilon \mathbf{L}_x(z))|\frac{J_{\mathbf{L}_x(z)}}{|\Gamma_x^1|}dz,$$
undoing again the change of variable, we obtain,
$$\left|\int_{B(0,1)}\nabla_{\mathbf{y}}w_\varepsilon(x, \varepsilon \mathbf{L}_x(z))\varepsilon\frac{\partial}{\partial x}(\mathbf{L}_x(z))\frac{J_{\mathbf{L}_x(z)}}{|\Gamma_x^1|}dz\right|\leq C\frac{\varepsilon}{|\Gamma_x^\varepsilon|}\int_{\Gamma_x^\varepsilon}|\nabla_{\mathbf{y}}w_\varepsilon(x, \mathbf{y})|d\mathbf{y}.$$
We estimate the last term as follows,
$$\int_{B(0,1)}w_\varepsilon(x, \varepsilon \mathbf{L}_x(z))\frac{\partial}{\partial x}\left(\frac{J_{\mathbf{L}_x(z)}}{|\Gamma_x^1|}\right)dz=$$
$$\int_{B(0,1)}\left(w_\varepsilon(x, \varepsilon \mathbf{L}_x(z))- (M_\varepsilon w_\varepsilon)(x)\right)\frac{\partial(J_{\mathbf{L}_x(z)}/|\Gamma_x^1|)}{\partial x}(z) dz +$$
$$+ (M_\varepsilon w_\varepsilon)(x)\int_{B(0,1)}\frac{\partial(J_{\mathbf{L}_x(z)}/|\Gamma_x^1|)}{\partial x}(z) dz. $$
Since,
$$\int_{B(0,1)}\frac{\partial (J_{\mathbf{L}_x(z)}/|\Gamma_x^1|)}{\partial x}(z) dz=\frac{d}{dx}\left(\frac{1}{|\Gamma_x^1|}\underbrace{\int_{B(0,1)}J_{\mathbf{L}_x(z)}(z)dz}_{|\Gamma_x^1|}\right)=0,$$
then, we have
$$\int_{B(0,1)}w_\varepsilon(x, \varepsilon \mathbf{L}_x(z))\frac{\partial}{\partial x}\left(\frac{J_{\mathbf{L}_x(z)}}{|\Gamma_x^1|}\right)dz=$$
$$=\int_{B(0,1)}\left(w_\varepsilon(x, \varepsilon \mathbf{L}_x(z))- (M w_\varepsilon)(x)\right)\frac{\partial(J_{\mathbf{L}_x(z)}/|\Gamma_x^1|)}{\partial x}(z) dz.$$
As before,  undoing the change of variable and taking account that $\left|\frac{\partial(J_{\mathbf{L}_x(z)}/|\Gamma_x^1|)}{\partial x}\right|\leq C$, we obtain
$$\left|\int_{B(0,1)}w_\varepsilon(x, \varepsilon \mathbf{L}_x(z))\frac{\partial}{\partial x}\left(\frac{J_{\mathbf{L}_x(z)}}{|\Gamma_x^1|}\right)dz\right|\leq C\frac{1}{|\Gamma_x^\varepsilon|}\int_{\Gamma_x^\varepsilon}|w_\varepsilon(x, \mathbf{y})- M_\varepsilon w_\varepsilon(x)|d\mathbf{y}.$$
Then, if we put together the three obtained estimates, we have

$$\left|M_\varepsilon\frac{\partial w_\varepsilon}{\partial x}-(M_\varepsilon w_\varepsilon)'\right|\leq C\frac{\varepsilon}{|\Gamma_x^\varepsilon|}\int_{\Gamma_x^\varepsilon}\nabla_{\mathbf{y}}w_\varepsilon(x, \mathbf{y})d\mathbf{y} + C\frac{1}{|\Gamma_x^\varepsilon|}\int_{\Gamma_x^\varepsilon}(w_\varepsilon(x, \mathbf{y})- M_\varepsilon w_\varepsilon(x))d\mathbf{y}.$$
So,
$$|\tilde{I}_1|=\left|\int_{Q_\varepsilon}\left(M_\varepsilon\frac{\partial w_\varepsilon}{\partial x}-\frac{\partial M_\varepsilon w_\varepsilon}{\partial x}\right)\frac{\partial v_\varepsilon}{\partial x}\right|\leq \int_0^1C\varepsilon\int_{\Gamma_x^\varepsilon}(\nabla_{\mathbf{y}}w_\varepsilon)v'_\varepsilon d\mathbf{y}dx+$$
$$+\int_0^1C\int_{\Gamma_x^\varepsilon}(w_\varepsilon- M_\varepsilon w_\varepsilon)v'_\varepsilon d\mathbf{y}dx.$$

Applying the H\"{o}lder inequality, $|\tilde{I}_1|$ can be estimated as follows,
$$\left|\int_{Q_\varepsilon}\left(M_\varepsilon\frac{\partial w_\varepsilon}{\partial x}-\frac{\partial M_\varepsilon w_\varepsilon}{\partial x}\right)\frac{\partial v_\varepsilon}{\partial x}\right|\leq C\varepsilon\|\nabla_{\mathbf{y}}w_\varepsilon\|_{L^2(Q_\varepsilon)}\|v'_\varepsilon\|_{L^2(Q_\varepsilon)}+ $$
$$+C\|w_\varepsilon-E_\varepsilon M_\varepsilon w_\varepsilon\|_{L^2(Q_\varepsilon)}\|v'_\varepsilon\|_{L^2(Q_\varepsilon)}.$$
By Lemma \ref{normaproyeccionextension}, 
$$\|w_\varepsilon-E_\varepsilon M_\varepsilon w_\varepsilon\|_{L^2(Q_\varepsilon)}\leq \sqrt{\beta}\varepsilon \|\nabla_{\mathbf{y}}w_\varepsilon\|_{L^2(Q_\varepsilon)},$$
so,
$$\left|\int_{Q_\varepsilon}\left(M_\varepsilon\frac{\partial w_\varepsilon}{\partial x}- (M_\varepsilon w_\varepsilon)'\right) v'_\varepsilon ds\right|\leq C \varepsilon\|v'_\varepsilon\|_{L^2(Q_\varepsilon)}\|\nabla_{\mathbf{y}}w_\varepsilon\|_{L^2(Q_\varepsilon)}.$$
To estimate the norm $\|v'_\varepsilon\|_{L^2(Q_\varepsilon)}$ we proceed as follows. We know that $v_\varepsilon$ is the solution of
\begin{equation}
\left\{
\begin{array}{r l r}
-\frac{1}{g}(g {v_\varepsilon}_x)_x + \mu v_\varepsilon\;&\; = M_\varepsilon f_\varepsilon, \;&\;\textrm{in}\quad (0, 1)\\
{v_\varepsilon}_x(0)\;&\;=0,\;&\;{v_\varepsilon}_x(1)=0.
\end{array}
\right.
\end{equation} 
Then, for $x\in(0,1)$, $v_\varepsilon$ satisfies,
$$-(gv'_\varepsilon)'+g\mu v_\varepsilon=gM_\varepsilon f_\varepsilon.$$
If we multiply by $v_\varepsilon$ and integrate by parts, we obtain,
$$\int_0^1g(v'_\varepsilon)^2dx+\mu\int_0^1gv_\varepsilon^2dx = \int_0^1(M_\varepsilon f_\varepsilon) gv_\varepsilon dx \leq$$
$$\stackrel{\textrm{H\"{o}lder ineq.}}{\leq}\left(\int_0^1(M_\varepsilon f_\varepsilon)^2dx\right)^{\frac{1}{2}}\left(\int_0^1(gv_\varepsilon)^2dx\right)^{\frac{1}{2}}\leq$$
$$\leq \frac{1}{4\delta}\int_0^1(M_\varepsilon f_\varepsilon)^2dx+\delta\int_0^1(gv_\varepsilon)^2dx.$$
Then,
$$\int_0^1g(v'_\varepsilon)^2dx +(\mu-\delta) \int_0^1gv_\varepsilon^2dx\leq\frac{1}{4\delta}\int_0^1(M_\varepsilon f_\varepsilon)^2dx\leq\|M_\varepsilon f_\varepsilon\|^2_{L^2(0,1)}$$
So,
$$\int_0^1g(x)(v'_\varepsilon)^2dx+ \int_0^1gv_\varepsilon^2dx\leq C\|M_\varepsilon f_\varepsilon\|^2_{L^2(0,1)}.$$
And,
$$\|v'_\varepsilon\|^2_{L^2(Q_\varepsilon)}=\int_{Q_\varepsilon}(v'_\varepsilon)^2ds=\int_0^1\int_{\Gamma_x^\varepsilon}(v'_\varepsilon)^2d\mathbf{y} dx=\varepsilon^{d-1}\int_0^1 g(x)(v'_\varepsilon)^2dx\leq$$
$$\leq \varepsilon^{d-1}C\|M_\varepsilon f_\varepsilon\|^2_{L^2(0,1)}.$$
Then,$$\left|\int_{Q_\varepsilon}\left(M_\varepsilon\frac{\partial w_\varepsilon}{\partial x}-(M_\varepsilon w_\varepsilon)'\right) v'_\varepsilon ds\right|\leq C \varepsilon\|v'_\varepsilon\|_{L^2(Q_\varepsilon)}\|\nabla_{\mathbf{y}}w_\varepsilon\|_{L^2(Q_\varepsilon)}\leq $$
$$\leq C\varepsilon\varepsilon^{\frac{d-1}{2}}\|M_\varepsilon f_\varepsilon\|_{L^2(0,1)}\|\nabla_{\mathbf{y}} w_\varepsilon\|_{L^2(Q_\varepsilon)}=$$
$$=C\varepsilon^{\frac{d+1}{2}}\|M_\varepsilon f_\varepsilon\|_{L^2(0,1)}\|\nabla_{\mathbf{y}}w_\varepsilon\|_{L^2(Q_\varepsilon)}\leq C\varepsilon^{d+1}\|M_\varepsilon f_\varepsilon\|^2_{L^2(0,1)} + \frac{1}{4}\|\nabla_{\mathbf{y}}w_\varepsilon\|^2_{L^2(Q_\varepsilon)}.$$
Note that,
$$\|\nabla_{\mathbf{y}}w_\varepsilon\|^2_{L^2(Q_\varepsilon)}\leq\|\nabla w_\varepsilon-\nabla v_\varepsilon\|^2_{L^2(Q_\varepsilon)},$$
so,
$$\left|\int_{Q_\varepsilon}\left(M_\varepsilon\frac{\partial w_\varepsilon}{\partial x}-(M_\varepsilon w_\varepsilon)'\right)v'_\varepsilon ds\right|\leq C\varepsilon^{d+1}\|M_\varepsilon f_\varepsilon\|^2_{L^2(0,1)}+\frac{1}{4}\|\nabla w_\varepsilon-\nabla v_\varepsilon\|^2_{L^2(Q_\varepsilon)}.$$

And $\tilde{I}_3$ can be estimated as follows,
$$\tilde{I}_3=\int_{Q_\varepsilon}(f_\varepsilon-M_\varepsilon f_\varepsilon)(w_\varepsilon-v_\varepsilon)ds=\int_{Q_\varepsilon}(f_\varepsilon-M_\varepsilon f_\varepsilon)(w_\varepsilon-M_\varepsilon w_\varepsilon)ds + $$
$$+\underbrace{\int_{Q_\varepsilon}(f_\varepsilon-M_\varepsilon f_\varepsilon)(M_\varepsilon w_\varepsilon-v_\varepsilon)ds}_{=0},$$
by the H\"{o}lder inequality,
$$|\tilde{I}_3|=\left|\int_{Q_\varepsilon}(f_\varepsilon-M_\varepsilon f_\varepsilon)(w_\varepsilon-M_\varepsilon w_\varepsilon)ds\right|\leq\|f_\varepsilon-M_\varepsilon f_\varepsilon\|_{L^2(Q_\varepsilon)}\|w_\varepsilon-M_\varepsilon w_\varepsilon\|_{L^2(Q_\varepsilon)}.$$
Again, by Lemma \ref{normaproyeccionextension},
$$\|w_\varepsilon-M_\varepsilon w_\varepsilon\|^2_{L^2(Q_\varepsilon)}\leq\beta\varepsilon^2\|\nabla_{\mathbf{y}} w_\varepsilon\|^2_{L^2(Q_\varepsilon)},$$
so,
$$|\tilde{I}_3|=\left|\int_{Q_\varepsilon}(f_\varepsilon-M_\varepsilon f_\varepsilon)(w_\varepsilon-M_\varepsilon w_\varepsilon)ds\right|\leq\|f_\varepsilon-M_\varepsilon f_\varepsilon\|_{L^2(Q_\varepsilon)}\sqrt{\beta}\varepsilon\|\nabla_{\mathbf{y}}w_\varepsilon\|_{L^2(Q_\varepsilon)}.$$
If we join all the estimates, then
$$\lambda_\varepsilon=\frac{1}{2}\int_{Q_\varepsilon}|\nabla w_\varepsilon-\nabla v_\varepsilon|^2 + |w_\varepsilon-v_\varepsilon|^2+\varepsilon^{d-1}\tau_\varepsilon+\theta_\varepsilon,$$
where,
$$|\theta_\varepsilon|=|\tilde{I}_1 - \tilde{I}_3|\leq C\varepsilon^{d+1}\|M_\varepsilon f_\varepsilon\|^2_{L^2(0,1)}+\frac{1}{4}\|\nabla w_\varepsilon-\nabla v_\varepsilon\|^2_{L^2(Q_\varepsilon)}+$$
$$+\|f_\varepsilon-M_\varepsilon f_\varepsilon\|_{L^2(Q_\varepsilon)}\sqrt{\beta}\varepsilon\|\nabla_{\mathbf{y}}w_\varepsilon\|_{L^2(Q_\varepsilon)}.$$
With this,
$$\lambda_\varepsilon\geq \frac{1}{2}\int_{Q_\varepsilon}(|\nabla w_\varepsilon-\nabla v_\varepsilon|^2 + |w_\varepsilon-v_\varepsilon|^2)ds+\varepsilon^{d-1}\tau_\varepsilon-C\varepsilon^{d+1}\|M_\varepsilon f_\varepsilon\|^2_{L^2(0,1)}$$
$$-\frac{1}{4}\|\nabla w_\varepsilon-\nabla v_\varepsilon\|^2_{L^2(Q_\varepsilon)}-\|f_\varepsilon-M_\varepsilon f_\varepsilon\|_{L^2(Q_\varepsilon)}\sqrt{\beta}\varepsilon\|\nabla_{\mathbf{y}}w_\varepsilon\|_{L^2(Q_\varepsilon)}.$$
By Lemma \ref{normaproyeccionextension} $\|M_\varepsilon f_\varepsilon\|_{L^2_g(0,1)}\leq \varepsilon^{\frac{1-d}{2}}\|f_\varepsilon\|_{L^2(Q_\varepsilon)}$, then
$$\lambda_\varepsilon\geq \frac{1}{4}\int_{Q_\varepsilon}(|\nabla w_\varepsilon-\nabla v_\varepsilon|^2 + |w_\varepsilon-v_\varepsilon|^2)ds+$$
$$+\varepsilon^{d-1}\tau_\varepsilon-C\varepsilon^{2}\|f_\varepsilon\|^2_{L^2(Q_\varepsilon)}-\|f_\varepsilon-M_\varepsilon f_\varepsilon\|_{L^2(Q_\varepsilon)}\sqrt{\beta}\varepsilon\|\nabla_{\mathbf{y}}w_\varepsilon\|_{L^2(Q_\varepsilon)}.$$
If we put everything together,
$$\varepsilon^{d-1}\tau_\varepsilon\geq\lambda_\varepsilon\geq\frac{1}{4}\int_{Q_\varepsilon}(|\nabla w_\varepsilon-\nabla v_\varepsilon|^2 + |w_\varepsilon-v_\varepsilon|^2)ds+\varepsilon^{d-1}\tau_\varepsilon-C\varepsilon^2\|f_\varepsilon\|^2_{L^2(Q_\varepsilon)}-$$
$$-\|f_\varepsilon-M_\varepsilon f_\varepsilon\|_{L^2(Q_\varepsilon)}\sqrt{\beta}\varepsilon\|\nabla_{\mathbf{y}}w_\varepsilon\|_{L^2(Q_\varepsilon)}\geq$$
$$\geq\frac{1}{4}\int_{Q_\varepsilon}(|\nabla w_\varepsilon-\nabla v_\varepsilon|^2 + |w_\varepsilon-v_\varepsilon|^2)ds+\varepsilon^{d-1}\tau_\varepsilon-C\varepsilon^2\|f_\varepsilon\|^2_{L^2(Q_\varepsilon)}-$$ 
$$\frac{\beta\varepsilon^2}{2} \|f_\varepsilon-M_\varepsilon f_\varepsilon\|^2_{L^2(Q_\varepsilon)}-\frac{1}{2}\|\nabla_{\mathbf{y}}w_\varepsilon\|^2_{L^2(Q_\varepsilon)},$$
so,
$$\|\nabla w_\varepsilon-\nabla v_\varepsilon\|^2_{L^2(Q_\varepsilon)}+\|w_\varepsilon-v_\varepsilon\|^2_{L^2(Q_\varepsilon)}\leq C\varepsilon^2\|f_\varepsilon\|^2_{L^2(Q_\varepsilon)}+\frac{\beta\varepsilon^2}{2} \|f_\varepsilon-M_\varepsilon f_\varepsilon\|^2_{L^2(Q_\varepsilon)},$$

and so,
$$\|w_\varepsilon-E_\varepsilon v_\varepsilon\|_{H^1(Q_\varepsilon)}\leq C\varepsilon\|f_\varepsilon\|_{L^2(Q_\varepsilon)},$$
that is,
\begin{equation}\label{estimacionresolvente}
\|u_\varepsilon-Ev_\varepsilon\|_{H^1_{\bm\varepsilon}(Q)}\leq C\varepsilon\|f_\varepsilon\|_{L^2(Q)}.\end{equation}
\begin{flushright}$\blacksquare$\end{flushright}

\section{Appendix: Shadowing and distance of attractors in $\mathbb{R}^m$}\label{shadowing}
In this section we introduce some concepts of the known Shadowing theory. The aim of this theory is to study the relationship between trajectories of a given dynamical system and trajectories of a perturbation of it. These techniques allow us to relate the distance of attractors in $\mathbb{R}^m$ with the distance of the corresponding time one maps for an ordinary differential equation and an appropriate perturbation of it. This result is described in Proposition \ref{LipschitzShadowingUniform}. Shadowing theory plays an important role in our work. Most of definitions we present below, can be found in \cite{Al-Nayef&Diamond&Kloeden&Co}. 

Throughout this section we will denote by $X$ a Banach space, $B\subset X$ a subset which may be bounded or unbounded and $T$ a nonlinear map, no necessary continuous or differentiable.

\begin{defi}\label{trajectory}
 A {\bf negative trajectory} of a map $T$ is a sequence  $\mathbf{x}_- = \{x_n\}_{n\in\mathbb{Z}^-}\subset B $ such that 
 $$x_{n+1}=T(x_n),\qquad\textrm{ for}\quad n\in\mathbb{Z}^-.$$
 \end{defi}
\begin{defi}\label{pseudotrajectory}
Let $\delta\geq 0$. A {\bf negative} $\bm\delta${\bf-pseudo-trajectory} of $T$ is a sequence $\mathbf{y} =\{y_n\}_{n\in\mathbb{Z}^-} \subset B$ with 
$$||y_{n+1} - T(y_n)|| \leq \delta, \quad \textrm{for}\quad n\in  \mathbb{Z}^- .$$
\end{defi}
We denote by $Tr^-(T, K,\delta)$ the set of all negative $\delta$-pseudo-trajectories of $T$ in $K \subset B$.  Note that 
a negative $0$-pseudo trajectory is a negative trajectory and that we always have the following inclusion
$$Tr^-(T, K, 0) \subset Tr^-(T, K, \delta).$$

An important class of negative $\delta$-pseudo-trajectories of a map $T$ are given by trajectories, $\{y_n\}_{n\in\mathbb{Z}^-}$, of maps $\varphi$, with $\varphi:B\to X$, such that for any $x\in B$, $\|T(x)-\varphi(x)\|\leq\delta$. This follows directly from the fact that
$$\|T(y_n)-y_{n+1}\| = \|T(y_n)-\varphi(y_n)\|\leq \delta\quad for \quad n\in\mathbb{Z}^-.$$
That is,
$$\bigcup_{\|T-\varphi\|\leq\delta}Tr^- (\varphi, B, 0)\subset Tr^-(T, B, \delta).$$

In this work we are going to need to compare the set of negative $\delta$-pseudo trajectories and the set of negative trajectories of a map $T:B\to B$.   An appropriate concept for this is the concept of ``Lipschitz Shadowing''. Hence, we consider the space, $l^p(X)$, for $1\leq p<\infty$, of all infinite negative sequences $\{x_n\}_{n\in\mathbb{Z}^-}$ such that
$$\|x\|_{l^p}=\Big(\sum_{j=1}^\infty|x_j|^p\Big)^{1/p}<\infty,$$
and $l^\infty(X)$ the Banach space given by the sequences $\mathbf{x_-}=\{x_n\}_{n\in\mathbb{Z}^-}$ with $x_n\in X$ and $\|x_n\|_X\leq C$ for all $n\in\mathbb{Z}^-$. That is,
$$l^\infty(X)=\{\mathbf{x_-}=\{x_n\}_{n\in\mathbb{Z}^-} : x_n\in X,\quad \|x_n\|_X\leq C\quad \forall n\in\mathbb{Z}^-\},$$
with $C>0$ a constant and the norm
$$\|\mathbf{x_-}\|_{l^\infty(X)}=\sup\{\|x_n\|_X: n\in\mathbb{Z}^-\}.$$
It is well known these spaces with these norms are Banach spaces.
\begin{defi}
 A negative sequence $\mathbf{x_-}=\{x_n\}_{n\in\mathbb{Z}^-}$ $\bm\varepsilon${\bf-shadows} a negative sequence $\mathbf{y_-}=\{y_n\}_{n\in\mathbb{Z}^-}$ if and only if, $$\|\mathbf{x_-}-\mathbf{y_-}\|_{l^\infty(X)}\leq\varepsilon.$$
So, this property is commutative, that is, $\mathbf{x_-}$ $\varepsilon$-shadows $\mathbf{y_-}$ if and only if $\mathbf{y_-}$ $\varepsilon$-shadows $\mathbf{x_-}$
\end{defi}
If for a given sequence $\mathbf{y_-}\in l^\infty(X)$ and $\eps>0$ we define $$B_\eps(\mathbf{y_-})=\{ \mathbf{x}_-=\{x_n\}_{n\in \mathbb{Z}^-}:  \|\mathbf{x_-} - \mathbf{y_-}\|_{l^\infty(X)}<\eps\},$$
then, we can write that a negative sequence $\mathbf{x_-}=\{x_n\}_{n\in\mathbb{Z}^-}$ $\bm\varepsilon${\bf-shadows} a sequence $\mathbf{y_-}=\{y_n\}_{n\in\mathbb{Z}^-}$ if $\mathbf{x_-}\in B_\eps(\mathbf{y_-})$. 
Finally, the main concept we want to present in this section is the following.
\begin{defi}\label{lipschitzshadowing}
The map $T$ has the {\bf Lipschitz Shadowing} property on $K\subset B$, if there exist constants $L, \delta_0> 0$ such that for any $0<\delta\leq \delta_0$ and any negative $\delta$-pseudo-trajectory of $T$ in $K$ is $(L\delta)$-shadowed by a negative trajectory of $T$ in $B$, that is, 
$$Tr^-(T,K,\delta)\subset B_{L\delta}(Tr^-(T,B,0)).$$
\end{defi}
 \bigskip
 
 All these concepts allow us to present the following result.
\begin{prop}\label{LipschitzShadowingUniform}

Let 
\begin{equation}\label{origineq}
\dot{x} = f(x), 
\end{equation}
be a dissipative Morse-Smale system. We perturbe it
\begin{equation}\label{pertureq}
\dot{x} = f_\varepsilon (x),
\end{equation}
with $\varepsilon\geq 0$ such that 
$$\bar{T}_\varepsilon\rightarrow \bar{T}_0\qquad as \quad \varepsilon\rightarrow 0,$$
in the $C^1$ topology and
$$\bar{T}_{\varepsilon}, \bar{T}_0: \mathbb{R}^m\rightarrow\mathbb{R}^m,$$
with $\bar{T}_0$ and $\bar{T}_\varepsilon$ the time one maps of the discrete dynamical systems generated by the evolution equations  (\ref{origineq}) and (\ref{pertureq}), respectively. Assume that for each $\varepsilon>0$, $\bar{T}_0$ and $\bar{T}_\varepsilon$ have global attractors $\bar{\mathcal{A}}_0$ and $\bar{\mathcal{A}}_\varepsilon$, respectively. Then we have

\begin{equation}\label{dist-attractors-Apendix}
dist_H (\bar{\mathcal{A}}_0, \bar{\mathcal{A}}_\varepsilon) \leq \bold{C}\|\bar{T}_0-\bar{T}_\varepsilon\|_\infty,
\end{equation}

with $\bold{C}$ independent of $\varepsilon$. 
\end{prop}
\par\medskip

\paragraph{\sl Proof.} 
Since $\bar{T}_0$ is a Morse-Smale map, in \cite{PilyuginShaDyn} the author proves that then, there exists a neighborhood of $\bar{T}_0$ in the $C^1$ topology, $\Theta$, and numbers $L, \delta_0$ such that, for any map $T'\in\Theta$, $T'$ has the Lipschitz Shadowing property on $\mathcal{N}(\bar{\A})$ with constants $\delta_0$, $L>0$.

On one side, since $\bar{T}_0$ has the Lipschitz Shadowing property on $\mathcal{N}(\bar{\A})$ with parameters $L, \delta_0$, then any negative $\delta$-pseudo-trajectory of $\bar{T}_0$ in $\mathcal{N}(\bar{\A})$, $\delta\leq\delta_0$, is $L\delta$-shadowed by a negative trajectory of $\bar{T}_0$ in $\mathbb{R}^m$, i.e.,
$$Tr^-(\bar{T}_0, \mathcal{N}(\bar{\A}),\delta)\subset B_{L\delta}(Tr^-(\bar{T}_0, \mathbb{R}^m, 0)).$$
Take $\varepsilon$ small enough such that $\|\bar{T}_\varepsilon-\bar{T}_0\|_\infty = \delta$ and $\bar{\mathcal{A}}_\varepsilon\subset \mathcal{N}(\bar{\A})$. We consider $r^\varepsilon\in\bar{\mathcal{A}}_\varepsilon,$ with 
$$\mathbf{r}_-^{\bm{\varepsilon}}=\{r^\varepsilon_n\}_{n\in\mathbb{Z}^-}=\{\bar{T}^n_\varepsilon(r^\varepsilon): n\in\mathbb{Z}^-\}\subset\bar{\mathcal{A}}_\varepsilon$$ its negative trajectory under the dynamical system generated by $\bar{T}_\varepsilon$,
$$\mathbf{r}_-^{\bm{\varepsilon}}\in {Tr}^-(\bar{T_\varepsilon}, \bar{\mathcal{A}}_\varepsilon,0).$$ As we have mentioned above, $\mathbf{r}_-^{\bm{\varepsilon}}$ is a negative $\delta$-pseudo-trajectory of $\bar{T}_0$ in $\bar{\mathcal{A}}_\varepsilon\subset\mathcal{N}(\bar{\A})$,
$$\mathbf{r}_-^{\bm{\varepsilon}}\in {Tr}^-(\bar{T}_0, \mathcal{N}(\bar{\A}),\delta).$$
So, there exist $\mathbf{r}_-=\{r_n\}_{n\in\mathbb{Z}^-}\in Tr^-(\bar{T}_0, \mathbb{R}^m, 0)$ such that,
$$\|r^\varepsilon_n-r_n\|\leq L\delta,$$
for all $n$ for which $\mathbf{r_-}$ is defined. Since
$$\|r_n\|\leq \|r^\varepsilon_n\|+ L\delta,$$
we conclude that $\mathbf{r_-}$ is bounded and for this reason $\mathbf{r_-}\in\bar{\A}$. With this
$$dist(r^\varepsilon, \bar{\A})\leq L\delta= L \|\bar{T}_\varepsilon-\bar{T}_0\|_\infty.$$
Since $r^\varepsilon\in\bar{\mathcal{A}}_\varepsilon$ has been chosen in an arbitrary way, we have
\begin{equation}\label{upperAuto}
dist(\bar{\mathcal{A}}_\varepsilon, \bar{\A})\leq L\delta = L \|\bar{T}_\varepsilon-\bar{T}_0\|_\infty,
\end{equation}
where $L$ is independent of $\varepsilon$.

On the other side, since any $T'\in\Theta$ has the Lipschitz Shadowing property on $\mathcal{N}(\bar{\A})$ of constants $L, \delta_0$, we take $\varepsilon>0$ small enough such that $\bar{T}_\varepsilon\in\Theta$ and 
$$\|\bar{T}_0-\bar{T}_\varepsilon\|_\infty\leq\delta<\delta_0.$$
With this, we take $r_0\in\bar{\mathcal{A}}_0$ and its negative trajectory under $\bar{T}_0$, $\mathbf{r_-}=\{r_n\}_{n\in\mathbb{Z}^-}=\{\bar{T}^n_0(r_0): n\in\mathbb{Z}^-\}$. As we have mentioned before, $\mathbf{r_-}$ is a negative $\delta$-pseudo-trajectory of $\bar{T}_\varepsilon$ in $\bar{\A}\subset\mathcal{N}(\bar{\A})$. Since we have chosen an small $\varepsilon$ such that $\bar{T}_\varepsilon\in\Theta$, then $\bar{T}_\varepsilon$ has the Lipschitz Shadowing property on $\mathcal{N}(\bar{\A})$ with parameters $L, \delta_0$, that is,
$$\mathbf{r_-}\in B_{L\|\bar{T}_0-\bar{T}_\varepsilon\|_\infty}({Tr}^-(\bar{T}_\varepsilon,\mathbb{R}^m,0)).$$
So, there exist $\mathbf{r}_-^{\bm{\varepsilon}}\in {Tr}^-(\bar{T}_\varepsilon, \mathbb{R}^m,0)$ such that
$$\|r_n-r^\varepsilon_n\|\leq L\|\bar{T}_0-\bar{T}_\varepsilon\|_\infty,$$
for all $n$ for which $\mathbf{r}_-^{\bm{\varepsilon}}$ is defined. Thus $\mathbf{r}_-^{\bm{\varepsilon}}$ is bounded. For that, $\mathbf{r}_-^{\bm{\varepsilon}}\in\bar{\mathcal{A}}_\varepsilon$ and also we have
$$\|r_0-r_0^\varepsilon\|\leq L\|\bar{T}_0-\bar{T}_\varepsilon\|_\infty,$$
with $r_0$ and $r_0^\varepsilon$ the $n=0$ elements of the sequences $\mathbf {r_-}$ and $\mathbf{r}_-^{\bm{\varepsilon}}$ respectively. That is
$$dist(r_0,\bar{\mathcal{A}}_\varepsilon)\leq L\|\bar{T}_0-\bar{T}_\varepsilon\|_\infty.$$
Finally, again since $r_0\in\bar{\A}$ have been chosen in an arbitrary way we conclude
\begin{equation}\label{lowerAuto}
dist(\bar{\A}, \bar{\mathcal{A}}_\varepsilon)\leq L\|\bar{T}_0-\bar{T}_\varepsilon\|_\infty.
\end{equation}

If we put together $(\ref{upperAuto})$ and $(\ref{lowerAuto})$ we obtain the desired estimate,
$$dist_H (\bar{\mathcal{A}}_0, \bar{\mathcal{A}}_\eps) \leq L\|\bar{T}_0-\bar{T}_\varepsilon\|_\infty,$$ 
with $L$ independent of $\varepsilon$.
\begin{flushright}$\blacksquare$\end{flushright}

\begin{re} 
Observe that the constant $C$ in \eqref{dist-attractors-Apendix} is the constant from the Lischitz Shadowing property of the map $T$.
\end{re}
An immediate consequence of the result above is the following 

\begin{cor}\label{MS-dist}
Let $T:\mathbb{R}^m\rightarrow\mathbb{R}^m$ be a {\bf Morse-Smale (gradient like)} map which has a global attractor $\mathcal{A}$. Then, there exists a neighborhood $\Theta$ of $T$ in the $C^1(\mathcal{N}(\mathcal{A}), \mathbb{R}^m)$ topology so that, for any $\,T_1, T_2\in\Theta$ with $\mathcal{A}_1$, $\mathcal{A}_2$ its respective attractors, we have
$$dist_H(\mathcal{A}_1, \mathcal{A}_2)\leq L\|T_1-T_2\|_{L^\infty(\mathcal{N}(\mathcal{A}), \mathbb{R}^m)},$$
with $L$ the Lipschitz Shadowing constant from the map $T$.
\end{cor}
\paragraph{\sl Proof.} As we mentioned in the proof of the proposition above, since $T$ is a Morse-Smale map, in \cite{PilyuginShaDyn} the author proves that there exists a neighborhood $\Theta$ of $T$ in the $C^1$ topology and numbers $L, \delta_0$ such that, for any map $T'\in\Theta$, $T'$ has the Lipschitz Shadowing property on $\mathcal{N}(\bar{\A})$ with constants $\delta_0$, $L>0$. The rest of the proof follows the same lines as the proof of the proposition above. 
\begin{flushright}$\blacksquare$\end{flushright}

\end{document}